\documentclass[11pt,reqno,bold]{paper}
\usepackage[letterpaper,left=25mm,right=25mm,top=25mm,bottom=25mm,headsep=5mm,footskip=7mm]{geometry}

\usepackage{amsmath}
\usepackage{amssymb}
\usepackage{amsfonts}
\usepackage{amsthm}
\usepackage{mathrsfs}

\usepackage{enumerate}
\usepackage{color}
\definecolor{darkblue}{rgb}{0.0,0.0,0.4}
\usepackage[pdfauthor={Tsogtgerel Gantumur},pdftitle={prescribed scalar curvature mean curvature problem},pdfsubject={Conformal geometry},colorlinks,citecolor=darkblue,linkcolor=darkblue,urlcolor=darkblue]{hyperref}

\usepackage[backend=bibtex,style=numeric,maxcitenames=6,maxbibnames=6,doi=false,url=false,isbn=false]{biblatex}
\DeclareNameAlias{sortname}{first-last}
\bibliography{prescribed}

\newtheorem{theorem}{Theorem}[section]
\newtheorem{proposition}[theorem]{Proposition}
\newtheorem{lemma}[theorem]{Lemma}
\newtheorem{corollary}[theorem]{Corollary}
\theoremstyle{remark}

\newcommand{\R}{{\mathbb{R}}}

\newcommand{\Tr}{{\gamma}}
\newcommand{\Yam}{\mathscr{Y}}

\title{A prescribed scalar and boundary mean curvature problem and the Yamabe classification on asymptotically Euclidean manifolds with inner boundary}
\author{Vladmir Sicca and Gantumur Tsogtgerel}
\institution{McGill University}
\date{\today}                                           

\begin{document}

\maketitle

\begin{abstract}
We consider the problem of finding a metric in a given conformal class with prescribed non-positive scalar curvature and non-positive boundary mean curvature
on an asymptotically Euclidean manifold with inner boundary. We obtain a necessary and sufficient condition in terms of a conformal invariant of the zero sets of the target curvatures for the existence of solutions to the problem and use this result to establish the Yamabe classification of metrics in those manifolds with respect to the solvability of the prescribed curvature problem.
\end{abstract}

\section{Introduction}
\label{2sec:main}

The classical Yamabe problem for closed Riemannian manifolds can be seen as the question of finding a metric in a given conformal class with constant scalar curvature, or prescribing constant scalar curvature in a given conformal class (a good review can be found in \cite{LP87}). The problem was studied for compact manifolds with boundary by Escobar, where he prescribes a constant scalar curvature in the interior of the manifold and a (different) constant mean curvature on its boundary (see for example \cite{Esco96a}). In the case of $(M,\partial M)$ a manifold with boundary, the problem reduces to finding solutions to the equation
\begin{equation}
    \begin{cases}
    -\Delta u+\frac{n-2}{4(n-1)}Ru=\frac{n-2}{4(n-1)}C u^{2\bar q-1},\ in\ M,\\
    \Tr \partial_\nu u+\frac{n-2}{2}H\Tr u=\frac{n-2}{2}D(\Tr u)^{\bar q},\ on\ \partial M,
\end{cases}
\end{equation}
with $\Delta$ the metric Laplacian that, in the euclidean case reduces to $\Delta u=\displaystyle\sum_i \partial^2_{i} u$, $R$ the original scalar curvature in the interior of $M$, $H$ the original mean curvature on the boundary $\partial M$ and $C$, $D$ the constants prescribed as the new scalar curvature and mean curvature on the boundary, respectively. In this case, the conformal factor corresponding to the new metric would be $u^{2\bar q-2}$, with $\bar q=\frac{n}{n-2}$. Then solutions can be found by minimizing the functional
\begin{equation}
    E(u)=\int_\Omega |\nabla u|dV_g+\frac{n-2}{4(n-1)}\int_\Omega Ru^2dV_g+\frac{n-2}{2}\int_\Sigma H(\Tr u)^2 d\sigma_g
\end{equation}
over an appropriate function space. For manifolds without boundary, both the equation and the functional to be considered are the same, except that we need to drop the boundary components.

An important result in the study of the aforementioned problems on compact manifolds is that there is a classification of metrics with respect to the sign of the minimizing energy. In the case of manifolds without boundary, one can only find solutions of constant curvature if the sign of the curvature agrees with the sign of the minimizer. In the case of manifolds with boundary, Escobar's results show that the picture is a little more complicated, although with a similar flavor. The full classification, called the Yamabe classification of the manifold, reads as follows:
\begin{itemize}
    \item If the minimizing energy is negative, it is possible to realize $C\equiv 0$ and $D<0$ or $C<0$ and $D$ a constant of any sign,
    \item If the minimizing energy is positive, it is possible to realize $C\equiv 0$ and $D>0$ or $C>0$ and $D$ a constant of any sign,
    \item If the minimizing energy is zero, only the case $C=D=0$ can be solved.
\end{itemize}

With this problem solved, the next step can be to address the situation where $C$ and $D$ are not constants, but general functions $R'$ and $H'$ that will represent the scalar curvature and the mean curvature, respectively, after the conformal transformation. In this case the conformal factor is not determined by the minimization of $E$, but the the Yamabe classification of the manifold still plays a role in determining the existence of solutions to the problem. In the compact case, \cite{DM15} for closed manifolds and \cite{ST21} for manifolds with boundary get existence results for prescribing non-positive scalar curvature and mean curvature on the boundary when the minimizing energy is negative. The paper \cite{DM15} also provides the complete Yamabe classification of asymptotically euclidean manifolds without boundary with respect to the existence of solutions to the prescribed curvature problem for non-positive scalar curvatures. 

In this paper we address the problem of prescribing scalar curvature and mean curvature on the boundary of an asymptotically euclidean manifold with (compact) boundary. To this end we adapt the techniques developed in \cite{DM15} for the case without boundary and do the necessary additions to bring into account the boundary terms. The root of the argument are variational techniques as the ones used to solve the problem in the compact case in \cite{ST21}, but with extra difficulties given the fact that the manifolds have infinite volume, so the usual Sobolev embedding theorems do not apply. Most of them are circumvented by the use of weighted Sobolev spaces as in \cite{Max05} and \cite{DM15}, only composed of functions that decay near infinity, that have similar embedding properties to Sobolev spaces in compact manifolds, as described in Section \ref{2sec:weighted_spaces}.

Our culminating results are similar to the case of manifolds without boundary in \cite{DM15}. Theorem \ref{2t:general_existence_result} gives a necessary and sufficient condition to the existence of solutions to the prescribed non-positive scalar curvature and non-positive mean curvature on the boundary as the combination of the zero sets of the target curvatures being Yamabe positive, in a sense defined in Section \ref{2sec:yamabe} for pairs of subsets of the interior of the manifold and of the boundary. This also leads to a classification of conformal classes of metrics on those manifolds similar to the one in the case without boundary seen in Theorem \ref{2t:yamabe_classification}, that can be summarized as:
\begin{itemize}
    \item The manifold is Yamabe positive if, and only if, any non-positive scalar curvature can be realized along with any non-positive mean curvature on the boundary;
    \item The manifold is Yamabe zero if, and only if, any non-positive scalar curvature can be realized along with any non-positive mean curvature on the boundary as long as we do not have both vanishing identically at the same time;
    \item The manifold is Yamabe negative otherwise.
\end{itemize}
The only remarkable difference to the case without boundary is for Yamabe zero manifolds, that can realize a vanishing scalar curvature in the manifold as long as the prescribed mean curvature on the boundary is not identically zero as well. This shows a consonance with the classification results of Escobar and the existence results discussed in \cite{ST21} for compact manifolds with boundary, where the conditions on the boundary presents less constraints to the existence results than the conditions in the interior of the manifold.

As a remark, results like these are expected to have consequences in the study of existence of solutions to the Lichnerowicz equation, that emerges when one searches for solutions to the Einstein constraint equations in a given conformal metric. Both in the case of asymptotically euclidean manifolds without boundary (in \cite{Max09}) and in the case of compact manifolds with boundary (in \cite{HT13}) there are results that reduce the existence of solutions to the Lichnerowicz equation to existence of solutions to the prescribed curvature problem in some cases. We do not know of similar results for asymptotically euclidean manifolds with boundary, but if those are found Theorem \ref{2t:general_existence_result} would automatically present a condition for the existence solution to the Lichnerowicz equation as well.

\subsection{Outline of the paper}

This paper is structured as follows. In Section \ref{2sec:weighted_spaces} we present our precise definition of an asymptotically euclidean manifold and state basic properties of weighted Sobolev spaces. In Section \ref{2sec:gluing} we present a structure for our asymptotically euclidean manifolds with boundary that splits it into a compact part and the ends of the manifold, and introduce a construction that will be useful in translating some estimates from the compact manifolds with boundary and from the asymptotically euclidean without boundary cases to our problem on asymptotically euclidean manifolds with boundary. In Section \ref{2sec:yamabe} we do the necessary calculations to build the relative Yamabe invariants on an asymptotically euclidean manifold with boundary and their relevant properties to prove the main results of the paper. In Section \ref{2sec:prescribed} we prove the existence results for the prescribed curvature problem and the Yamabe classification of asymptotically euclidean manifolds with boundary with respect to the possibility of prescription of nonpositive curvatures.
\section{Weighted Sobolev Spaces}
\label{2sec:weighted_spaces}

Throughout this paper, let $(M,\partial M, g)$ be a smooth, connected, $n$-dimensional manifold with boundary $\partial M$ that can be decomposed as $M\cup\partial M=K\cup E_1\cup ...\cup E_p$, with $K$ a compact set, $\partial M\subset K$, and each of the $E_i$'s diffeomorphic to $\mathbb{R}^n\setminus B_1(0)$ with $E_i\cap E_j=\emptyset$ for $i\neq j$. Let us also specify a diffeomorphism $f_i:\mathbb{R}^n\setminus B_1(0)\to E_i$ for each $i$. In this case we say each $E_i$ is an {\it{end}} of the manifold and that $M$ has $p$ ends. We want $M$ to be asymptotically euclidean, that is, if $\eta$ is the usual euclidean metric in $\mathbb{R}^n$, then $f_i^\ast g(x)\to \eta(x)$ as $\|x\|\to\infty$, where $f_i^{\ast}g$ is the pullback of the metric $g$ through the diffeomorphism $f_i$. We need to qualify this convergence, though, and for that we will use weighted Sobolev spaces.

Let us fix a smooth function $\rho\geq 1$ over $M$ such that, for each $i$, $f_i^\ast \rho=\|x\|$ outside of a compact set in $\mathbb{R}^n$ and a smooth Riemannian metric $\hat g$ over $M$ such that $f_i^\ast\hat g=\eta$ also outside a compact set in $\mathbb{R}^n$. We define
\begin{equation}
\label{2eq:w_k_p_delta_norm}
    \|u\|^p_{W^{k,p}_\delta(M)}:=\sum_{j=0}^{k}\left\vert\left\vert\rho^{-\delta-\frac{n}{p}+j}\left\vert\hat\nabla^j u\right\vert_{\hat g}\right\vert\right\vert^p_{L^p(M,\hat g)}
\end{equation}
for $k$ a non-negative integer and for a number $\delta\in\mathbb{R}$. In the definition, $\hat\nabla^j$ is the $j$-th covariant derivative with respect to the Levi-Civita connection induced by $\hat g$ and $u\in L^1_{loc}(M,g)$. As usual, we say $u\in W^{k,p}_{\delta}(M)$ if, and only if, $\|u\|_{W^{k,p}_\delta(M)}<+\infty$ and we denote $W^{0,p}_\delta(M):=L^p_\delta(M)$. We will say that $(M,\partial M, g)$ is asymptotically euclidean and $g$ is a $W^{k,p}_\delta$ asymptotically euclidean metric if $g-\hat g\in W^{k,p}_\delta(M)$. As mentioned in \cite{DM15}, if $g\in W^{k,p}_\delta(M)$ and $0\leq k\leq 2$, the norm one obtains by replacing $\hat g$ by $g$ in (\ref{2eq:w_k_p_delta_norm}) is equivalent to the one defined using $\hat g$, so under this restriction throughout the text we will be using the one that is more convenient to the situation discussed without further explanation.

Some properties of weighted spaces should be highlighted because they will be used throughout the paper. First, straight from the definitions, we have that $L^q(M)=L^q_{-\frac{n}{q}}(M)$. In particular, if we define 
\begin{equation}
\delta^\ast=\frac{2-n}{2},
\end{equation}
we have that $L^{2\bar q}=L^{2\bar q}_{\delta^\ast}$ and this index behaves as a critical index for the weighted spaces in some sense and thus will be an important threshold in this paper. One of the reasons is that
\begin{equation}
    \|u\|^2_{W^{1,2}_{\delta^\ast}(M)}=\|\nabla u\|^2_{L^2(M)}+\|u\|^2_{L^2_{\delta^\ast}(M)},
\end{equation}
so if $u\in W^{1,2}_{\delta^\ast}(M)$, then $\|\nabla u\|^2_{L^2(M)}<+\infty$ and $\delta^\ast$ is the largest index guaranteeing that, that is, it is not necessarily true for $u\in W^{1,2}_\delta(M)$. For example, if $\delta>\delta^\ast$, $M=\mathbb{R}^n\setminus B_1(0)$ and $u=r^\alpha$ for suitably chosen exponents $\alpha$. Of course, as $\rho\geq 1$, $\delta<\delta^\ast$ implies $W^{1,2}_\delta(M)\subset W^{1,2}_{\delta^\ast}(M)$.

Another important property of weighted spaces is that $L^p_\delta(M)\hookrightarrow L^q_{\delta'}(M)$ if $p\geq q$ and $\delta<\delta'$. In terms of the critical index, an important consequence is that $L^{2\bar q}(M)=L^{2\bar q}_{\delta^\ast}(M)\hookrightarrow L^q_{\delta}(M)$ for $q\leq 2\bar q$ and $\delta>\delta^\ast$. To finish the list of important embeddings that mimic the case for compact manifolds, item 3 of Lemma 1 in \cite{Max05} implies $W^{1,2}_{\delta^\ast}(M)\hookrightarrow L^{2\bar q}(M)$ and the usual Sobolev embeddings can be used with attention to some restrictions on $\delta$.

Finally, two important inequalities are proven in Lemma 2.1 of \cite{DM15} for functions in $W^{1,2}_{\delta^\ast}(M)$, that we state here. For $u\in W^{1,2}_{\delta^\ast}(M)$, there are $C_1$, $C_2$ positive constants such that
\begin{equation}
\label{2eq:poincare_inequality}
    \|u\|_{L^2_{\delta^\ast}(M)}\leq C_1 \|\nabla u\|_{L^2(M)}
\end{equation}
and
\begin{equation}
\label{2eq:sobolev_inequality}
    \|u\|_{L^{2\bar q}(M)}\leq C_2\|\nabla u\|_{L^2(M)}.
\end{equation}

\section{Gluing Framework}
\label{2sec:gluing}

In this paper we try to adapt the results proven in \cite{ST21} for compact manifolds with boundary and in \cite{DM15} for asymptotically euclidean manifolds without boundary to the case of asymptotically euclidean manifolds {\it{with}} boundary. One might expect that nothing really new is necessary to allow the results to be transferred from the previous cases to the new one, since the compact $K$ can be taken as a compact manifold with boundary and the new phenomena in the asymptotically euclidean case should be controlled by the behavior of functions in the ends (near infinity), not affected by the eventual presence of boundaries in the compact part of the manifold. That is in fact true in some situations, in particular when establishing a priori estimates for the quantities we will be studying later, although it is harder to do when we deal with the solutions to differential equations.

Given the many cases this intuition can be used, in this section we develop a construction that will be evoked in the proofs of some results throughout the paper when we can decompose the estimates into their counterparts -- on one hand in the compact manifold with boundary case and, on the other, in the asymptotically euclidean without boundary case. The notation developed will also be useful in other cases when we consider either only the compact part or the end part of a function, even if we do not use the whole decomposition argument. Hence, it will be used freely throughout the paper.

Let $\tilde K\subset M\cup\partial M$ be a relatively open set such that $K\subset \tilde K$ and the closure of $\tilde K$ is a smooth compact manifold with boundary. We also call $E:=E^\circ_1\cup...\cup E^{\circ}_p$. Notice that $\tilde K\cup E=M\cup\partial M$ and we can pick $\{\xi_E, \xi_{\tilde K}\}$ a smooth partition of unity subordinated to $E\cup \tilde K$ in the usual way.

We also define $(\tilde M_p, \tilde g)$ as a smooth Riemannian manifold without boundary obtained by isometrically attaching $p$ copies of $\mathbb{R}^n\setminus B_1(0)$ to an n-sphere with $p$ holes in a way that $\tilde g$ is smooth. If we call $\tilde f_i:\mathbb{R}^n\setminus B_1(0)\to \tilde M$ the distinguished  diffeomorphisms associated to the isometric gluing (that is, $\tilde f_i$ is an isometry), we see that $(\tilde M_p, \tilde g)$ is $W^{k,p}_\delta$ asymptotically euclidean for any $k$, $p$, $\delta$ as $\tilde g$ plays the role of both the usual metric in the manifold and of $\hat g$ in the previous section. Accordingly, we call $\tilde E=\bigcup \left(\tilde f_i\circ f^{-1}_i\right)(E_i)$. We finish our setting by defining $\tilde \rho\leq 1$ a smooth function over $\tilde M$ such that $\tilde f_i^\ast\tilde\rho=f_i^\ast\rho$ for each $i$.

Now, if we have a function $u$ over $M$ with support in $E$, we can define a new function $\tilde u$ over $\tilde M$ by $\tilde u(x)=u\left(\left(f_i \circ\tilde f_i^{-1}\right)(x)\right)$ if $x\in \left(\tilde f_i\circ f_i^{-1}\right)(E_i)$ and $\tilde u(x)=0$ otherwise. Then we have that, since $u$ has support in $E$
\begin{equation}
    \int_M \rho^{-\delta-\frac{n}{p}+j}|\hat\nabla^j u|dV_{\hat g}=\sum_{i=1}^p\int_{f_i^{-1}(E_i)}f_i^\ast\rho^{-\delta-\frac{n}{p}+j}|\partial_j f_i^\ast u|\sqrt{\det \hat g_{ij}}d\mathbf{x}
\end{equation}
with the coordinate quantities corresponding to the charts defined by the $f_i$'s. But we have $\det \hat g=1$ except on a compact set. So there are $C_1,\ C_2$ such that
\begin{equation}
    C_1\sum_{i=1}^p\int_{f_i^{-1}(E_i)}f_i^\ast\rho^{-\delta-\frac{n}{p}+j}|\partial_j f_i^\ast u|d\mathbf{x}\leq \int_M \rho^{-\delta-\frac{n}{p}+j}|\hat\nabla^j u|dV_{\hat g}\leq C_2\sum_{i=1}^p\int_{f_i^{-1}(E_i)}f_i^\ast\rho^{-\delta-\frac{n}{p}+j}|\partial_j f_i^\ast u|d\mathbf{x}
\end{equation}
and, by the way we defined $\tilde \rho$ and $\tilde u$, and since the $\tilde f_i$'s are isometries
\begin{equation}
\label{2eq:equivalence_metrics}
    C_1\|\tilde u\|^p_{W^{s,p}_\delta(\tilde M)}\leq \| u\|^p_{W^{s,p}_\delta( M)}\leq C_2\|\tilde u\|^p_{W^{s,p}_\delta(\tilde M)}
\end{equation}
as long as $u\in W^{s,p}_\delta(M)$. As a notation, quantities satisfying the property illustrated in equation (\ref{2eq:equivalence_metrics}) will be denoted as $\|\tilde u\|^p_{W^{s,p}_\delta(\tilde M)}\sim \| u\|^p_{W^{s,p}_\delta( M)}$.
\section{The relative Yamabe invariant}
\label{2sec:yamabe}

From now on, we assume that $(M,\partial M, g)$ is a $W^{2,p}_\tau$ asymptotically euclidean manifold with boundary of dimension $n\geq 3$, $p>\frac{n}{2}$ and $\tau<0$. Denote its scalar curvature by $R\in L^p_{\tau-2}(M)$ and the mean extrinsic curvature at the boundary by $H\in W^{1-\frac{1}{p},p}(\partial M)$, taking the mean curvature with respect to the inner normal. Let $\Omega\subset M$ and $\Sigma \subset \partial M$ be relatively measurable sets and consider the functional $E_g:W^{1,2}_{\delta^\ast}(M)\to \mathbb{R}$ defined by

\begin{equation}
E_g(u)=\int_\Omega |\nabla u|dV_g+\frac{n-2}{4(n-1)}\int_\Omega Ru^2dV_g+\frac{n-2}{2}\int_\Sigma H(\Tr u)^2 d\sigma_g,    
\end{equation}
with $\Tr:W^{1,2}_\delta(M)\to W^{\frac{1}{2}, 2}(\partial M)$ the trace map and $dV_g$ and $d\sigma_g$ the volume forms induced by $g$ on $M$ and $\partial M$ respectively. We will omit the subscript when the metric in use is understood.

We will look for minimizers for $E$ in certain function spaces. Specifically, let
\begin{equation}
    W^{1,2}_{\delta^\ast}(\Omega, \Sigma):=\overline{\left\{u\in C^1_{\delta^\ast} (M)\cap C^0(\bar M):u|_{M\setminus \Omega}\equiv 0, \Tr(u)|_{\partial M\setminus \Sigma}\equiv 0\right\}}
\end{equation}

with the closure taken over $W^{1,2}_{\delta^\ast}(M)$ and define
\begin{equation}
    B^{q,r}_{b}(\Omega,\Sigma)=\left\{u\in W^{1,2}_{\delta^\ast} (\Omega,\Sigma): \|u\|^q_{L^q(\Omega)}+b\|\Tr u\|^r_{L^r(\Sigma)}=1 \right\}.
\end{equation}

In this section we will study the quantities

\begin{equation}
    \Yam^{q,r}_{b}(\Omega,\Sigma):=\inf_{u\in B^{q,r}_{b}(\Omega,\Sigma)}E(u).
\end{equation}

\begin{lemma}
\label{2l:bounds_on_curvature_terms}
Let $\delta>\delta^\ast$ and let $u\in W^{1,2}_{\delta^\ast}(M)$, $\Omega$ be a subset of $M$ and $\Sigma$ be a subset of $\partial M$. Then, for any $\epsilon>0$, there is a constant $K_\epsilon>0$ such that
\begin{equation}
\label{2eq:control_R_term}
    \left\vert\int_\Omega Ru^2 dV_g\right\vert\leq \epsilon\|\nabla u\|^2_{L^2(\Omega)}+K_\epsilon\|u\|^2_{L^2_\delta(\Omega)},
\end{equation}
and
\begin{equation}
\label{2eq:control_H_term}
    \left\vert\int_\Sigma H(\Tr u)^2 d\sigma_g\right\vert\leq \epsilon\|\Tr u\|^2_{W^{\frac{1}{2}, 2}(\Sigma)}+K_\epsilon\|\Tr u\|^2_{L^2(\Sigma)}.
\end{equation}
\begin{proof}
Using the construction from section \ref{2sec:gluing} we have:
\begin{equation}
\label{2e:Ru2_estimate_beginning}
\begin{array}{ccl}
    \left\vert \int_\Omega Ru^2dV_g\right\vert & = &\left\vert \int_\Omega R((\xi_{\tilde K}+\xi_E)u)^2dV_g\right\vert \\
     & \leq & \left\vert \int_\Omega R\xi_{\tilde K}^2u^2dV_g\right\vert+\left\vert 2\int_\Omega R\xi_E\xi_{\tilde K}u^2dV_g\right\vert+\left\vert \int_\Omega R\xi_E^2u^2dV_g\right\vert.
\end{array}
\end{equation}
We are going to look for bounds for each of the terms.

For the first term, using equation (14) in \cite{ST21} for $\tilde K$ we get
\begin{equation*}
    \begin{array}{ccl}
        \left\vert \int_\Omega R\xi_{\tilde K}^2u^2dV_g\right\vert & =&\left\vert \int_{\Omega\cap \tilde K} R\xi_{\tilde K}^2u^2dV_g\right\vert \\
        &\leq & \epsilon \|\nabla(\xi_{\tilde K}u)\|^2_{L^2(\Omega\cap\tilde K)}+C_\epsilon\|\xi_{\tilde K}u\|^2_{L^2(\Omega\cap \tilde K)} \\
        &\leq& \epsilon\|\xi_{\tilde K}\nabla u\|^2_{L^2(\Omega\cap\tilde K)}+\epsilon\|u\nabla\xi_{\tilde K}\|^2_{L^2(\Omega\cap\tilde K)}+C_\epsilon\|\xi_{\tilde K}u\|^2_{L^2(\Omega\cap \tilde K)}
    \end{array}
\end{equation*}
and since $\xi_{\tilde K}$ is a smooth partition of unity with $\xi_{\tilde K}\leq 1$, and $\tilde K$ is precompact:
\begin{equation}
\label{2e:Ru2_estimate:R_estimate_compact_part}
    \left\vert \int_\Omega R\xi_{\tilde K}^2u^2dV_g\right\vert\leq \epsilon\left(\|\nabla u\|^2_{L^2(\Omega\cap \tilde K)}+\max_{\tilde K}|\nabla\xi_{\tilde K}|^2\|u\|^2_{L^2(\Omega\cap \tilde K)}\right)+C_\epsilon\|u\|^2_{L^2(\Omega\cap \tilde K)}.
\end{equation}
For the second term we use the fact that the same estimate (14) in \cite{ST21} does not depend on the fact that $R$ is the scalar curvature on the manifold, but only on its regularity. So we can define $\tilde R:=2R\xi_E\xi_{\tilde K}$ and get
\begin{equation}
\label{2e:Ru2_estimate:mixed_terms_estimate}
    \left\vert\int_\Omega 2R\xi_E\xi_{\tilde K}u^2dV_g\right\vert=\left\vert\int_\Omega \tilde R u^2dV_g\right\vert\leq \epsilon \|\nabla u\|^2_{L^2(\Omega\cap \tilde K)}+C'_{\epsilon}\|u\|^2_{L^2(\Omega\cap\tilde K)}.
\end{equation}
For the third term, notice that $R\xi_E^2u^2$ is supported in $E$, so
\begin{equation}
    \left\vert \int_\Omega R\xi_E^2u^2dV_g\right\vert=\left\vert \int_{\Omega\cap E} R\xi_E^2u^2dV_g\right\vert\sim\left\vert \int_{\tilde M} \left(R\xi_E^2u^2\right)^{\sim}dV_{\tilde g}\right\vert.
\end{equation}
Since $\tilde M$ is an asymptotically euclidean manifold without boundary, we can use Lemma 3.1 from \cite{DM15} on $\tilde M$ to have
\begin{equation*}
\begin{array}{ccl}
    \left\vert \int_{\tilde M}\left(R\xi_E^2u^2\right)^{\sim}dV_{\tilde g}\right\vert & \leq & \epsilon \|\nabla (\xi_E u)^{\sim}\|^2_{L^2(\tilde M)}+C''_{\epsilon}\|(\xi_E u)^\sim\|^2_{L^2_\delta(\tilde M)}\\
    & = &\epsilon \|\nabla(\xi_Eu)^\sim\|^2_{L^2(\tilde E)}+C''_\epsilon\|(\xi_Eu)^\sim\|^2_{L^2_\delta(\tilde E)},
\end{array}
\end{equation*}
because $\xi_E$ is supported in $E$. Now, since $\xi_E$ is a smooth partition of unity and $\nabla \xi_E\equiv 0$ outside of the precompact set $E\cap \tilde K$, we have that there is $C_1$ such that
\begin{equation}
    \epsilon C_1\|\nabla(\xi_Eu)^\sim\|^2_{L^2(\tilde E)}+C''_\epsilon\|(\xi_Eu)^\sim\|^2_{L^2_\delta(\tilde E)}\leq \epsilon\left(\max_{\tilde K}|\nabla\xi_E|^2\|u\|^2_{L^2(\Omega\cap \tilde K)}+\|\nabla u\|^2_{L^2(E)}\right)+C''_\epsilon C_1\|u\|^2_{L^2_\delta(E)}
\end{equation}
and coming back to $M$ though the diffeomorphisms between the copies of $E$ in $M$ and $\tilde M$ we conclude finally that there is $C_2>0$ such that
\begin{equation}
\label{2e:Ru2_estimate:R_estimate_ends}
    \left\vert \int_\Omega R\xi_E^2u^2dV_g\right\vert\leq \epsilon C_2\left(\max_{\tilde K}|\nabla\xi_E|^2\|u\|^2_{L^2(\Omega\cap \tilde K)}+\|\nabla u\|^2_{L^2(\Omega)}\right)+C''_\epsilon C_2\|u\|^2_{L^2_\delta(\Omega)}.
\end{equation}
Then we can gather the results of (\ref{2e:Ru2_estimate:R_estimate_compact_part}), (\ref{2e:Ru2_estimate:mixed_terms_estimate}) and (\ref{2e:Ru2_estimate:R_estimate_ends}) into (\ref{2e:Ru2_estimate_beginning}):
\begin{equation}
    \left\vert \int_\Omega Ru^2dV_g\right\vert\leq (2+C_2)\epsilon\|\nabla u\|^2_{L^2(\Omega)}+\left((1+C_2)\epsilon\max_{\tilde K}|\nabla\xi_E|^2+C_\epsilon+C'_\epsilon\right)\|u\|^2_{L^2(\Omega\cap \tilde K)}+C''_\epsilon C_2\|u\|^2_{L^2_\delta(\Omega)}
\end{equation}
and since $\rho\geq 1$
\begin{equation}
    \left\vert \int_\Omega Ru^2dV_g\right\vert\leq 3\epsilon\|\nabla u\|^2_{L^2(\Omega)}+\left(2\epsilon\max_{\tilde K}|\nabla\xi_E|^2+C_\epsilon+C'_\epsilon+C''_\epsilon\right)\|u\|^2_{L^2_\delta(\Omega)}
\end{equation}
which can be relabeled into (\ref{2eq:control_R_term}).

To prove (\ref{2eq:control_H_term}), notice that $\xi_{\tilde K}u$ is supported in $\tilde K$, $\Tr (\xi_{\tilde K}u)|_{(\partial \tilde K\setminus \partial M)}\equiv 0$ and, since in a neighborhood of $\partial M$ we have $\xi_{\tilde K}u=u$, $\Tr (\xi_{\tilde K}u)|_{\partial M}=\Tr u$. So using equation (15) of \cite{ST21}
\begin{equation}
    \left\vert\int_{\partial \tilde K} H(\Tr (\xi_{\tilde K}u))^2 d\sigma_g\right\vert\leq \epsilon\|\Tr (\xi_{\tilde K}u)\|^2_{W^{\frac{1}{2}, 2}(\partial \tilde K)}+K_\epsilon\|\Tr (\xi_{\tilde K}u)\|^2_{L^2(\partial \tilde K)}
\end{equation}
which translates immediately to (\ref{2eq:control_H_term}).
\end{proof}
\end{lemma}

\begin{lemma}
\label{2l:control_L2_norms}
Let $u\in B^{q,r}_{b}(\Omega,\Sigma)$, with $2\leq q\leq 2\bar q$ and $2\leq r\leq\bar q+1$, $q>r$ and $\delta>\delta^\ast$. Then for any $\epsilon>0$ there is $C_\epsilon>0$ independent of $u$ such that
\begin{equation}
\label{2eq:bound_L2d_norm_by_gradient}
    \|u\|^2_{L^2_\delta(\Omega)}\leq \epsilon\|\nabla u\|^2_{L^2(\Omega)}+C_\epsilon,
\end{equation}
and
\begin{equation}
\label{2eq:bound_boundary_normy_by_gradient}
    \|\gamma u\|^2_{L^2(\Sigma)}\leq \epsilon \|\nabla u\|^2_{L^2(\Omega)}+C_\epsilon.
\end{equation}
\begin{proof}
First we show that our hypotheses allow one to prove (\ref{2eq:bound_boundary_normy_by_gradient}) by applying Lemma 2.2 of \cite{ST21} to $\xi_{\tilde K}u$. Notice that, since $u$ is supported in $\Omega\cup\Sigma$ and $\xi_{\tilde K}$ is supported in $\tilde K$
$$
\begin{array}{ccl}
    \|\xi_{\tilde K}u\|^q_{L^q(\Omega\cap\tilde K)}+b\|\gamma(\xi_{\tilde K}u)\|^r_{L^r(\partial \tilde K)} &=& \|\xi_{\tilde K}u\|^q_{L^q(\Omega\cap\tilde K)}+b\|\gamma u\|^r_{L^r(\Sigma)} \\
    &\leq & \left(\displaystyle\max_{\tilde K}|\xi_{\tilde K}|^q\right)\|u\|^q_{L^q(\Omega\cap\tilde K)}+b\|\gamma u\|^r_{L^r(\Sigma)} \\
    &=&\|u\|^q_{L^q(\Omega\cap\tilde K)}+b\|\gamma u\|^r_{L^r(\Sigma)}=1,
\end{array}
$$
so $\xi_{\tilde K}u$ satisfies the hypothesis of the proof of Lemma 2.2 in \cite{ST21}. Hence, given $\epsilon>0$, there is ${C_\epsilon>0}$ such that
$$
\begin{array}{ccl}
    \|\gamma u\|^2_{L^2(\Sigma)}=\|\gamma(\xi_{\tilde K}u)\|^2_{L^2(\partial \tilde K)} & \leq &\epsilon\|\nabla(\xi_{\tilde K}u)\|^2_{L^2(\Omega\cap \tilde K)}+C_\epsilon \\
     &\leq &\epsilon\left(\|u\nabla\xi_{\tilde K}\|^2_{L^2(\Omega\cap \tilde K)}+\|\xi_{\tilde K}\nabla u\|^2_{L^2(\Omega\cap\tilde K)}\right)+C_{\epsilon}\\
     &\leq &\epsilon\left(\displaystyle\max_{\tilde K}|\nabla\xi_{\tilde K}|^2\|u\|^2_{L^2(\Omega\cap \tilde K)}+\displaystyle\max_{\tilde K}|\xi_{\tilde K}|^2\|\nabla u\|^2_{L^2(\Omega\cap\tilde K)}\right)+C_{\epsilon}\\
     &\leq&\epsilon\left(D\displaystyle\max_{\tilde K}|\nabla\xi_{\tilde K}|^2+1\right)\|\nabla u \|^2_{L^2(\Omega\cap \tilde K)}+C_{\epsilon},
\end{array}
$$
with $D$ the constant associated to the Sobolev inequality. This proves (\ref{2eq:bound_boundary_normy_by_gradient}) up to a relabelling since $\Omega\cap \tilde K\subset \Omega$.

On the other hand, by H\"older's inequality
\begin{equation}
    \|u\|^2_{L^2_\delta(\Omega)}=\int_\Omega u^2\rho^{-2\delta-n}dV_{g}\leq \left[\int_{\Omega}|u|^q dV_{ g}\right]^{\frac{2}{q}}\left[\int_{\Omega}\rho^{-\frac{2\delta q}{q-2}-\frac{nq}{q-2}}dV_{ g}\right]^{\frac{q-2}{q}}
\end{equation}
Notice that the exponent of $\rho$ inside the integral is less than $-n$, so the second integral is finite. Let us call $M=\int_\Omega \rho^{-\frac{2\delta q+nq}{q-2}}dV_{ g}$. Since
\begin{equation}
    \|u\|^q_{L^q(\Omega)}+b\|\gamma u\|^r_{L^r(\Sigma)}=1<2,
\end{equation}
replacing in the inequality above we have
$$
\begin{array}{ccl}
    \|u\|^2_{L^2_\delta(\Omega)} & \leq & M^{\frac{q-2}{q}}\left[2-b\|\gamma u\|^r_{L^r(\Sigma)}\right]^{\frac{2}{q}}\\
     &\leq & M^{\frac{q-2}{q}}\left[2+|b|\|\gamma u\|^r_{L^r(\Sigma)}\right]^{\frac{2}{q}}\\
     &\leq& 2^{\frac{2}{q}}M^{\frac{q-2}{q}}+|b|^{\frac{2}{q}}M^{\frac{q-2}{q}}\|\gamma u\|^{\frac{2r}{q}}_{L^r(\Sigma)},
\end{array}
$$
because $\frac{2}{q}\leq 1$.

Next, by continuity of the trace operator, there is $D_1>0$ such that
\begin{equation}
    \|\gamma u\|^2_{L^r(\Sigma)}\leq D_1\|\xi_{\tilde K}u\|^2_{W^{1,2}(\Omega\cap \tilde K)}\leq D_2\|u\|^{2}_{W^{1,2}(\Omega\cap \tilde K)},
\end{equation}
for some constant $D_2$ as in the calculations used in the proof of (\ref{2eq:bound_boundary_normy_by_gradient}). So we can replace in the previous expression to get
\begin{equation}
\label{2e:L2d_norm_estimate_with_h1_norm_and_r_over_q_exponents}
    \|u\|^2_{L^2_\delta(\Omega)}\leq D_3+D_4\|u\|^{\frac{2r}{q}}_{W^{1,2}(\Omega\cap \tilde K)}\leq D_3+D_4\left(\left(\int_{\Omega\cap\tilde K}|\nabla u|^2dV_g\right)^{\frac{r}{q}}+\left(\int_{\Omega\cap\tilde K}u^2dV_g\right)^{\frac{r}{q}}\right)
\end{equation}
where $D_3=2^{\frac{2}{q}}M^{\frac{q-2}{q}}$ and $D_4=|b|^{\frac{2}{q}}M^{\frac{q-2}{q}}D_2^{\frac{r}{q}}$, the second inequality because $\frac{r}{q}\leq 1$.

Now, given $\epsilon_1>0$, there is $D_5(\epsilon)>0$ such that, if $t\geq 0$:
\begin{equation}
    t^{\frac{r}{q}}\leq \epsilon_1 t+D_5
\end{equation}
taking both $t=\int_{\Omega\cap\tilde K}|\nabla u|^2dV_g$ and $t=\int_{\Omega\cap\tilde K}u^2dV_g$ we can replace the respective terms in (\ref{2e:L2d_norm_estimate_with_h1_norm_and_r_over_q_exponents}) to get
\begin{equation}
    \|u\|^2_{L^2_\delta(\Omega)}\leq D_6+D_4\epsilon_1\left(\int_{\Omega\cap\tilde K}|\nabla u|^2dV_g+\int_{\Omega\cap\tilde K}u^2dV_g\right)
\end{equation}
with $D_6=D_3+2D_4D_5$. But $\tilde K$ is precompact and $\rho\geq 1$ is smooth, so
\begin{equation}
\begin{array}{ccl}
    \|u\|^2_{L^2_\delta(\Omega)} & \leq & D_6+D_4\epsilon_1\left(\int_{\Omega\cap\tilde K}|\nabla u|^2dV_g+\frac{1}{\min_{\tilde K}\rho^{-2\delta-n}}\int_{\Omega\cap \tilde K}u^2\rho^{-2\delta-n}dV_g\right)\\
     & \leq & D_6+D_4\epsilon_1\left(\int_{\Omega\cap\tilde K}|\nabla u|^2dV_g+\frac{1}{\min_{\tilde K}\rho^{-2\delta-n}}\|u\|^2_{L^2_{\delta}(\Omega)}\right).
\end{array}
\end{equation}
So, choosing $\epsilon_1>0$ such that $D_5:=\frac{D_4\epsilon_1}{\min_{\tilde K}\rho^{-2\delta-n}}<1$, we finally get
\begin{equation}
    \|u\|^2_{L^2_\delta(\Omega)}
      \leq  \frac{D_6}{1-D_5}+\frac{D_4\epsilon_1}{1-D_5}\left(\int_{\Omega\cap\tilde K}|\nabla u|^2dV_g+\frac{1}{\min_{\tilde K}\rho^{-2\delta-n}}\|u\|^2_{L^2_{\delta}(\Omega)}\right),
\end{equation}
which can be relabelled into the result.
\end{proof}
\end{lemma}

The following results will use the next lemma, which comes up often in the study of manifolds with boundary by allowing a normalization of functions in $W^{1,2}_{\delta^\ast}(\Omega,\Sigma)$ into functions in $B^{q,r}_{b}(\Omega,\Sigma)$.

\begin{lemma}
\label{2l:polynomial}
Let $q>r>1$, $a>0$ and $b$ be constants (if $b>-a$ we can have $q=r$), and let
\begin{equation}
f_b(x) = ax^q+bx^r ,
\end{equation}
where $b\in\R$ is a parameter.
Then we have the following.
\begin{enumerate}[(a)]
\item
The equation $f_b(x)=1$ has a unique positive solution $x_b>0$.
\item
The correspondence $b\mapsto x_{b}$ is continuous.
\item $b\mapsto x_b$ is a non-increasing function.
\end{enumerate}
\begin{proof}
The results can be easily seen by looking at the graph of $f_b(x)$. For $(c)$, just notice that $ax^q_b+bx^r_b=1$ implies $ax^q_b+b'x^r_b\leq 1$, so $x_{b'}\geq x_b$.
\end{proof}
\end{lemma}

\begin{lemma}
\label{2l:energy_control_gradient}
Under the hypothesis of Lemma \ref{2l:bounds_on_curvature_terms}, there are constants $K_1$, $K_2$ and $K_3$ such that
\begin{equation}
\label{2eq:gradient_control_by_energy_and_norms}
    \|\nabla u\|^2_{L^2(\Omega)}\leq K_1E(u)+K_2\|u\|^2_{L^2_\delta(\Omega)}+K_3\|\Tr u\|^2_{L^2(\Sigma)},
\end{equation}
and if $u\in B^{q,r}_b(\Omega,\Sigma)$, $2\leq q\leq2\bar q$, $2\leq r\leq\bar q+1$ and $q>r$, there are $C$, $K$ such that
\begin{equation}
\label{2eq:gradient_control_by_only_energy}
    \|\nabla u\|^2_{L^2(\Omega)}\leq CE(u)+K.
\end{equation}
As a consequence, $\Yam^{q,r}_{b}(\Omega, \Sigma)$ is finite unless $W^{1,2}_{\delta^\ast}(\Omega,\Sigma)=\{0\}$.
\begin{proof}
From (\ref{2eq:control_R_term}) we have that, given $\epsilon_1>0$, there is $K_{\epsilon_1}>0$ such that
\begin{equation}
\label{2e:bounding_Yam_R_term_control}
    \int_{\Omega} Ru^2dV_g\geq-\epsilon_1\|\nabla u\|^2_{L^2(\Omega)}-K_{\epsilon_1}\|u\|^2_{L^2_{\delta}(\Omega)}.
\end{equation}
Similarly, by (\ref{2eq:control_H_term}), if $\epsilon_2>0$, there is $K_{\epsilon_2}$ satisfying
\begin{equation}
    \int_{\Sigma}H(\Tr u)^2d\sigma_g\geq -\epsilon_2\|\Tr u\|^2_{W^{\frac{1}{2},2}(\Sigma)}-K_{\epsilon_2}\|\Tr u\|^2_{L^2(\Sigma)},
\end{equation}
which implies, if $C_1$ is the constant associated to the trace inequality in $\tilde K$,
\begin{equation}
    \int_{\Sigma}H(\Tr u)^2d\sigma_g\geq -\epsilon_2C_1\|\xi_{\tilde K} u\|^2_{L^{2}(\Omega\cap \tilde K)}-\epsilon_2C_1\|\nabla(\xi_{\tilde K}u)\|^2_{L^2(\Omega\cap \tilde K)}-K_{\epsilon_2}\|\Tr u\|^2_{L^2(\Sigma)},
\end{equation}
and since $\xi_{\tilde K}\leq 1$ is smooth and $\tilde K$ is compact
\begin{equation}
\label{2e:bounding_Yam_H_term_control}
    \int_{\Sigma}H(\Tr u)^2d\sigma_g\geq -\epsilon_2C_1\left(1+\max_{\tilde K}|\nabla \xi_{\tilde K}|^2\right)\|u\|^2_{L^{2}(\Omega\cap \tilde K)}-\epsilon_2C_1\|\nabla(u)\|^2_{L^2(\Omega\cap \tilde K)}-K_{\epsilon_2}\|\Tr u\|^2_{L^2(\Sigma)}.
\end{equation}
Now, we also have that
\begin{equation}
    \|u\|^2_{L^2_\delta(\Omega)}\geq\int_{\Omega\cap \tilde K}\left(\rho^{-\delta-\frac{n}{2}}u\right)^2dV_g\geq \min_{\tilde K}\left(\rho^{-2\delta-n}\right)\|u\|^2_{L^2(\Omega\cap \tilde K)}.
\end{equation}
So, replacing in (\ref{2e:bounding_Yam_H_term_control}) we have
\begin{equation}
    \int_{\Sigma}H(\Tr u)^2d\sigma_g\geq -\epsilon_2C_2\|u\|^2_{L^{2}_\delta(\Omega)}-\epsilon_2C_1\|\nabla(u)\|^2_{L^2(\Omega\cap \tilde K)}-K_{\epsilon_2}\|\Tr u\|^2_{L^2(\Sigma)}.
\end{equation}
with $C_2=\left(\displaystyle\min_{\tilde K}\left(\rho^{-2\delta-n}\right)\right)^{-1}C_1\left(1+\max_{\tilde K}|\nabla \xi_{\tilde K}|^2\right)$. Adding it to (\ref{2e:bounding_Yam_R_term_control}), we can choose $\epsilon_1$, $\epsilon_2$ small enough such that $K_1^{-1}:=1-\frac{n-2}{4(n-1)}\epsilon_1-\frac{n-2}{2}\epsilon_2C_2>0$ and have finally
\begin{equation}
    E(u)\geq K_1^{-1}\|\nabla u\|^2_{L^2(\Omega)}-K_2K_1^{-1}\|u\|^2_{L^2_\delta(\Omega)}-K_3K^{-1}_1\|\Tr u\|^2_{L^2(\Sigma)}
\end{equation}
with $K_2:=\left(\frac{n-2}{4(n-1)}K_{\epsilon_1}+\frac{n-2}{2}\epsilon_2C_2\right)K_1$ and $K_3:=\frac{n-2}{2}K_{\epsilon_2}K_1$, and that can be reorganized as (\ref{2eq:gradient_control_by_energy_and_norms}).

Now, if $u\in B^{q,r}_b(\Omega,\Sigma)$, using Lemma \ref{2l:control_L2_norms}, we have also that, for a given $\epsilon>0$:
\begin{equation}
    E(u)\geq \left(K_1^{-1}-K_2K_1^{-1}\epsilon-K_3K^{-1}_1\epsilon\right)\|\nabla u\|^2_{L^2(\Omega)}-(K_2K_1^{-1}+K_3K^{-1}_1)C_\epsilon
\end{equation}
which, by choosing $\epsilon$ small enough, can be reorganized as (\ref{2eq:gradient_control_by_only_energy}).

Besides, the result shows $E(u)$ is bounded from below. So $\Yam_{b}^{q,r}(\Omega,\Sigma)$ is finite unless $B_{b}^{q,r}(\Omega,\Sigma)$ is empty, which, by Lemma \ref{2l:polynomial}, happens if, and only if, $W^{1,2}_{\delta^\ast}(\Omega,\Sigma)$ is trivial.
\end{proof}
\end{lemma}

We also need a technical lemma to guarantee that this invariant is meaningful in the context of conformal transformations.

\begin{lemma}
The quantities $\Yam^{2\bar q, r}_{0}(\Omega,\Sigma)$ and $\Yam^{2\bar q, \bar q+1}_{b}(\Omega,\Sigma)$ are invariant under conformal transformations of the metric.
\begin{proof}
Assume $g'=\phi^{\frac{4}{n-2}}g$ is a $W^{k,p}_{\delta}$ metric, $\phi-1\in W^{k,p}_\delta(M)$. Then, if we have $u\in B^{2\bar q,r}_{0}(\Omega,\Sigma)$ with respect to $g'$, $\phi u\in W^{1,2}_{\delta^\ast}(\Omega,\Sigma)$ with respect to $g$ (see the proof of Lemma 3.7 in \cite{DM15} and item 4 of Lemma 1 in \cite{Max05}), $\phi u\in B^{2\bar q, r}_{0}(\Omega,\Sigma)$ with respect to $g$ and $E_{g'}(u)=E_g(\phi u)$ (see the proof of Lemma 2.3 in \cite{ST21}) and the invariance of $\Yam^{2\bar q, r}_{0}(\Omega,\Sigma)$ follows. The proof for $\Yam^{2\bar q, \bar q+1}_{b}(\Omega,\Sigma)$ is similar.
\end{proof}
\end{lemma}

To be able to meaningfully talk about the conformal invariance we have to establish under which circumstances $\Yam^{q,r}_{b}(\Omega,\Sigma)$ is independent of the indices. In the next sequence of results we will prove that the sign of $\Yam^{q,r}_{b}(\Omega,\Sigma)$ does not in fact depend on $b$ or $r$, so we will end up fixing $q=2\bar q$ and talk about pairs of sets $(\Omega,\Sigma)$ that are Yamabe positive, negative or zero and write $\Yam(\Omega,\Sigma)>0$, $<0$ or $=0$ respectively. The proofs hold for compact manifolds with boundary as well, and are sometimes simpler than the ones presented in \cite{ST21}, but in that case we could guarantee independence of sign with respect to $q$ using the usual embeddings of $L^p$ spaces on sets of finite measure that are not available for non-compact manifolds with boundary.

First, we deal with the case $\Yam(\Omega,\Sigma)<0$, which is easier.

\begin{proposition}
\label{2p:independence_in_the_negative_case}
If $q'$, $r'$ and $b'$ satisfy the conditions of Lemma \ref{2l:polynomial}, $\Yam^{q,r}_{b}(\Omega,\Sigma)<0$ implies $\Yam^{q',r'}_{b'}(\Omega,\Sigma)<0$ independently of $q$, $r$ and $b$.
\begin{proof}
If $\Yam^{q,r}_{b}<0$, there is $u\in W^{1,2}_{\delta^\ast}(\Omega,\Sigma)$ such that $E(u)<0$. Lemma \ref{2l:polynomial}, then, provides a function $ku\in B^{q',r'}_{b'}(\Omega,\Sigma)$ for some $k>0$ and $E(ku)=k^2E(u)<0$.
\end{proof}
\end{proposition}

Now we deal with $b$, and we have a nice result on the dependence of $\Yam^{q,r}_{b}(\Omega,\Sigma)$ on $b$.

\begin{lemma}
If $q\in [2,2\bar q]$, $r\in [2,\bar q+1]$, $q>r$, the map $\Yam^{q,r}:b\mapsto \Yam^{q,r}_{b}(\Omega,\Sigma)$ is continuous for $b\in\mathbb{R}$. If we accept $q=r$ we have to restrict the domain to $b\geq 0$.
\begin{proof}
In this proof we will omit the sets $(\Omega,\Sigma)$ to simplify the notation.

Assume $\{b_n\}_n$ is a sequence in the domain of $\Yam^{q,r}$, $b_n\to b$. Fixing $u\in B^{q,r}_{b}$, for each $n$ there is $k_{n}(u)>0$ such that
\begin{equation}
    \|k_nu\|^q_{L^q(\Omega)}+b\|\gamma(k_n u)\|^r_{L^r(\Sigma)}=k_n^q\|u\|^q_{L^q(\Omega)}+k_n^rb\|u\|^r_{L^r(\Sigma)}=1.
\end{equation}
Now, if we consider the functions $f_u:y\mapsto y\|\gamma u\|^r_{L^r(\Sigma)}$, which is continuous, and $g_u:d\mapsto x_d$ such that $\|u\|^q_{L^q(\Omega)}x_d^q+dx_d^r=1$, which is continuous by Lemma \ref{2l:polynomial}, $g_u\circ f_u$ is continuous and
\begin{equation}
    \left(g_u\circ f_u\right)(b_n)=k_n.
\end{equation}
So $k_n\to (g_u\circ f_u)(b)=1$, since $u\in B^{q,r}_{b}$. As a consequence, $E(k_nu)=k_n^2E(u)\to E(u)$.

Our goal now is to prove that $\Yam^{q,r}_{b_n}\to \Yam^{q,r}_{b}$.

First, notice that $E(k_n u)\geq \Yam^{q,r}_{b_n}$ which means that, since $E(k_n u)\to E(u)$, if $n$ is large, then $\Yam^{q,r}_{b_n}<E(u)+1$ and $\{\Yam^{q,r}_{b_n}\}_n$ is bounded from above. Also, since this is true for any $u\in B^{q,r}_{b}$, if there is a converging subsequence $\Yam^{q,r}_{b_p}\to L$, $\Yam^{q,r}_{b}\geq L$.

In fact, if $\Yam^{q,r}_{b_p}\to L$, $L=\Yam^{q,r}_{b}$. Assume it is not true, so there is $\epsilon>0$ such that $\Yam^{q,r}_{b}=L+\epsilon$. As a consequence, for $p$ sufficiently large, $\Yam^{q,r}_{b_p}<L+\frac{\epsilon}{3}$ and there is $u_p\in B^{q,r}_{b_p}$ a function satisfying $E(u_p)\leq L+\frac{\epsilon}{3}$.

But again, under the hypothesis of the theorem, there is $c_p=\left(g_{u_p}\circ f_{u_p}\right)(b)>0$ such that $c_pu_p\in B^{q,r}_{b}$. It follows that
\begin{equation}
    E(c_pu_p)=c_p^2E(u_p)\geq L+\epsilon,
\end{equation}
and so 
\begin{equation}
    L+\frac{\epsilon}{3}\geq E(u_p)\geq \frac{L+\epsilon}{c_p^2}.
\end{equation}
Assuming $L\geq 0$, rearranging the two extremes of the inequality we have
\begin{equation}
    c_p^2\geq \frac{3L+3\epsilon}{3L+\epsilon}=1+\frac{2\epsilon}{3L+\epsilon}.
\end{equation}
But from the definition of the $c_p$'s
\begin{equation}
    c_p^q\|u_p\|^q_{L^q(\Omega)}+c_p^rb\|\Tr u_p\|^r_{L^r(\Sigma)}=1=\|u_p\|^q_{L^q(\Omega)}+b_p\|u_r\|^r_{L^r(\Sigma)}
\end{equation}
which results in
\begin{equation}
    (c_p^q-1)\|u_p\|^{q}_{L^q(\Omega)}+(c_p^rb-b_p)\|\Tr u_p\|^r_{L^r(\Sigma)}=0
\end{equation}
and, since $c_p^q-1>0$, for any $p$
\begin{equation}
    0\geq c^r_pb-b_p>\left(1+\frac{2\epsilon}{3L+\epsilon}\right)^{\frac{2}{r}}b-b_p,
\end{equation}
which is absurd because $b_p\to b$. So $L=\Yam^{q,r}_{b}$. If $L< 0$ the argument is similar.

If we prove that $\{\Yam_{b_n}^{q,r}\}_n$ is bounded from below, the result follows from the fact that all subsequences will have a subsequence that converges, hence converges to $\Yam^{q,r}_{b}$. But, if $\Yam^{q,r}_{b_n}$ has no lower bound, choosing $L<\Yam^{q,r}_{b}$ there is a subsequence $\{\Yam^{q,r}_{b_p}\}_p$ satisfying that for each $p$ there is $u_p$ such that $u_p\in B^{q,r}_{b_p}$ and $E(u_p)\leq L+\frac{\epsilon}{3}$ and we can repeat the previous argument verbatim. So our sequence is bounded and hence $\Yam^{q,r}$ is continuous.
\end{proof}
\end{lemma}

\begin{lemma}
\label{2l:monotonicity_with_respect_to_b}
Under the hypothesis of the previous lemma, if $\Yam^{q,r}_{b}(\Omega,\Sigma)<0$ for some set of indices, the map $\Yam^{q,r}$ is non-decreasing. Otherwise, it is non-increasing.
\begin{proof}
Throughout the proof, $f_u$, $g_u$ will be as in the previous lemma.

Assume $u\in B^{q,r}_{b}(\Omega,\Sigma)$. Then, if $b'\geq b$, $f_u(b')\geq f_u(b)$. So, as proved in Lemma \ref{2l:polynomial}, we have that $g_u(f_u(b'))\leq g_u(f_u(b))$. But $u\in B^{q,r}_{b}(\Omega,\Sigma)$ implies $g_u(f_u(b))=1$, so $0<k_u:=g_u(f_u(b'))\leq 1$. In addition, we have
\begin{equation}
    E(k_uu)=k_u^2E(u).
\end{equation}

If $\Yam^{q,r}_{b}(\Omega,\Sigma)<0$, we saw in Proposition \ref{2p:independence_in_the_negative_case} that the same is true for any set of indices $q'$, $r'$ and $b'$, so we can look at what happens only with functions of negative energy and assume $E(u)<0$. Then $E(k_uu)\geq E(u)$. Since any function in $B^{q,r}_{b'}(\Omega,\Sigma)$ can be written as $k_uu$ for some $u\in B^{q,r}_{b}(\Omega,\Sigma)$, $\Yam^{q,r}_{b'}(\Omega,\Sigma)\geq \Yam^{q,r}_{b}(\Omega,\Sigma)$.

On the other hand, if $\Yam^{q,r}_{b}(\Omega,\Sigma)\geq 0$, $E(u)\geq 0$ for any $u\in W^{1,2}_{\delta^\ast}(\Omega,\Sigma)$, so in the situation above, $E(k_uu)\leq E(u)$. As $u$ is arbitrary, $\Yam^{q,r}_{b'}(\Omega,\Sigma)\leq \Yam^{q,r}_{b}(\Omega,\Sigma)$.
\end{proof}
\end{lemma}

With these results, we can complete the proof of independence of the sign of $\Yam^{q,r}_{b}$ with respect to $b$ with the following two results.

\begin{lemma}
With the usual hypothesis on the indices, if $\Yam^{q,r}_{0}(\Omega,\Sigma)=0$, then $\Yam^{q,r}_{b}(\Omega,\Sigma)=0$ for any $b\in\mathbb{R}$.
\begin{proof}
We know already that, under the hypothesis, $\Yam^{q,r}_{b}\geq0$. If $b>0$ the result follows from Lemma \ref{2l:monotonicity_with_respect_to_b}. Assume then $b<0$.

Since $\Yam^{q,r}_{b}(\Omega,\Sigma)=0$, there is a minimizing sequence $\{u_n\}_n\subset B^{q,r}_{0}(\Omega,\Sigma)$ such that $E(u_n)\to 0$ and $\|u_n\|^q_{L^q(\Omega)}=1$ for all $n$. On the other hand, for each $n$ there is $k_n>0$ a number such that $k_nu_n\in B^{q,r}_{b}(\Omega, \Sigma)$, that is
\begin{equation}
\label{2e:k_n_relation_negative_b_0_case}
    k_n^q+bk_n^r\|\Tr u_n\|^r_{L^r(\Sigma)}=1.
\end{equation}
Besides, $E(k_nu_n)=k_n^2E(u_n)$. So either $\{k_n\}_n$ is bounded, and $E(k_nu_n)\to 0$, finishing the proof, or $k_n\to +\infty$. But since we have $b<0$, $q>r$ and equation (\ref{2e:k_n_relation_negative_b_0_case}) implies $\|\Tr u_n\|_{L^r(\Sigma)}\to\infty$.

On the other hand, if $E(u_n)\to 0$, the sequence of the $\{E(u_n)\}_n$ is bounded and Lemma \ref{2l:energy_control_gradient} implies a bound on the sequence $\{\|\nabla u_n\|_{L^2(\Omega)}\}_n$, which implies a bound in $\{\|u_n\|_{L^2_\delta(\Omega)}\}_n$ because of Lemma \ref{2l:control_L2_norms}. So $\{\|\xi_{\tilde K}u_n\|_{W^{1,2}(\Omega)}\}_n$ is also a bounded sequence and, by the trace inequality, there is a uniform bound on the $\|\Tr u_n\|_{L^r(\Sigma)}$. Contradiction. So $\Yam^{q,r}_{b}(\Omega,\Sigma)=0$.
\end{proof}
\end{lemma}

The next lemma gives the reverse direction.

\begin{lemma}
Under the same hypothesis on the indices, $\Yam^{q,r}_{b}(\Omega,\Sigma)=0$ implies $\Yam^{q,r}_{0}(\Omega,\Sigma)=0$.
\begin{proof}
Again, monotonicity gives us the result in the case $b<0$. So assume $b>0$ and let's take another value $b'>0$. Let $\{u_n\}_n\subset B^{q,r}_{b}(\Omega,\Sigma)$ be such that $E(u_n)\to 0$ and $k_n$ be satisfy
\begin{equation}
    k^q_n\|u_n\|^q_{L^q(\Omega)}+k_n^rb'\|\Tr u_n\|^r_{L^r(\Sigma)}=1.
\end{equation}
Again, either $k_n\to+\infty$ or $\Yam^{q,r}_{b'}(\Omega,\Sigma)=0$. But since $b'>0$, if $k_n\to+\infty$, $\|u_n\|_{L^q(\Omega)}\to 0$ and $\|\Tr u_n\|_{L^r(\Sigma)}\to 0$, absurd because $\{u_n\}_n\subset B^{q,r}_{b,\delta}(\Omega,\Sigma)$. So $\Yam^{q,r}_{b'}=0$ for all $b'>0$ and, by continuity, $\Yam^{q,r}_{0}(\Omega,\Sigma)=0$.
\end{proof}
\end{lemma}

We have as a corollary that the sign of the relative Yamabe invariant does not depend on $b$.

\begin{proposition}
\label{2p:independence_wrt_b}
Under the usual assumptions on the indices, the sign of $\Yam^{q,r}_{b}(\Omega,\Sigma)$ does not depend on $b$.
\begin{proof}
The previous two lemmata guarantee the proposition is true in the zero case, while Lemma \ref{2p:independence_in_the_negative_case} guarantees the result holds in the negative case. The positive case follows by exclusion.
\end{proof}
\end{proposition}

We have immediately also independence with respect to $r$.

\begin{proposition}
\label{2p:independence_wrt_r}
Under the same conditions on the indices, the sign of $\Yam^{q,r}_{b}(\Omega,\Sigma)$ does not depend on $r$.
\begin{proof}
The sign of $\Yam^{q,r}_{b}(\Omega,\Sigma)$ is the same as the sign of $\Yam^{q,r}_{0}(\Omega,\Sigma)$. But the calculation of $\Yam^{q,r}_{0}(\Omega,\Sigma)$ does not depend on $r$.
\end{proof}
\end{proposition}

We finally define the {\textit{relative Yamabe invariant}} of the pair $(\Omega,\Sigma)$ as $\Yam^{2\bar q,\bar q+1}_1(\Omega,\Sigma)$.

Although the Yamabe invariant is naturally related to the classical Yamabe problem, it is not the easiest invariant to be dealt with in the following proofs and calculations. Instead, we will work with the relative, weighted version of the first eigenvalue of the Laplacian defined as
\begin{equation}
    \lambda_{\delta}(\Omega,\Sigma)=\inf_{u\in W^{1,2}_{\delta^\ast}(\Omega,\Sigma)\setminus\{0\}}\frac{E(u)}{\|u\|^2_{L^2_\delta(\Omega)}+\|\Tr u\|^2_{L^2(\Sigma)}}.
\end{equation}
The bridge between both invariants is given by proving that the sign of $\lambda_\delta(\Omega)$ is the same as the sign of $\Yam_1^{2\bar q, \bar q+1}$ for any $\delta>\delta^\ast$. But first we can prove that there are minimizers for $\lambda_\delta(\Omega)$.

\begin{proposition}
\label{2p:minimizer_lambda}
If $W^{1,2}_{\delta^\ast}(\Omega,\Sigma)\neq\{0\}$, there is $u\in W^{1,2}_{\delta^\ast}(\Omega, \Sigma)$, $u> 0$, such that 
\begin{equation}
    \lambda_\delta(\Omega,\Sigma)=\frac{E(u)}{\|u\|^2_{L^2_\delta(\Omega)}+\|\Tr u\|^2_{L^2(\Sigma)}}.
\end{equation}
\begin{proof}
Let $\{u_k\}_k\subset W^{1,2}_{\delta^\ast}(\Omega,\Sigma)$ be a minimizing sequence such that $\|u_k\|^2_{L^2_\delta(\Omega)}+\|\Tr u_k\|^2_{L^2(\Sigma)}=1$ for any $k$. Since $E(v)=E(|v|)$ for all $v$, we can assume $u_k\geq 0$ for all $k$. Then $\|u_k\|^2_{L^2_\delta(\Omega)}, \|\Tr u_k\|^2_{L^2(\Sigma)}\leq 1$ for any $k$ and as the sequence is minimizing, it is bounded in $W^{1,2}_{\delta^\ast}(\Omega, \Sigma)$ due to estimate (\ref{2eq:gradient_control_by_energy_and_norms}).

As a consequence, there is $u\in W^{1,2}_{\delta^\ast}(\Omega,\Sigma)$ such that 

i. $u_k\rightharpoonup u$ in $W^{1,2}_{\delta^\ast}(\Omega)$ and then,

ii. $u_k\to u$ in $L^2_{\delta}(\Omega)$ and,

iii. $\Tr u_k\to \Tr u$ in $L^2(\Sigma)$.

As a first consequence of ii and iii, $\|u\|^2_{L^2_\delta(\Omega)}+\|\Tr u\|^2_{L^2(\Sigma)}=1$. Also, i implies 
\begin{equation}
    \|u\|_{W^{1,2}_{\delta^\ast}(\Omega)}\leq\liminf \|u_k\|_{W^{1,2}(\Omega)},
\end{equation} 
and that combined with ii gives $\|\nabla u\|_{L^2(\Omega)}\leq \liminf \|\nabla u_k\|_{L^2(\Omega)}$.

On the other hand, as proved in \cite{ST21}, the map $v\mapsto\int_{\partial M}H(\Tr v)^2d\sigma_g$ is continuous in $L^2(\Sigma)$, so
\begin{equation}
    \int_{\partial M}H(\Tr u_k)^2d\sigma_g\to \int_{\partial M}H(\Tr u)^2d\sigma_g.
\end{equation}
Also, the proof of Lemma 3.1 in \cite{DM15} holds for our case and because of ii we have that
\begin{equation}
    \int_{M}Ru_k^2dV_g\to \int_{M}Ru^2dV_g.
\end{equation}
As a result, $E(u)\leq \liminf E(u_k)=\lambda_\delta(\Omega, \Sigma)$ and $u$ is the minimizer we are looking for. Also because of ii, $u_k\to u$ pointwise, so $u\geq 0$.

Moreover, if $u$ is such minimizer, it is a weak solution to the equation
\begin{equation}
\label{2eq:linear_problem}
    \begin{cases}
    -\Delta u+\frac{n-2}{4(n-1)}Ru=\lambda u,\ in\ \Omega,\\
    \Tr\partial_\nu u+\frac{n-2}{2}H\Tr u=\lambda\Tr u,\ in\Sigma, 
    \end{cases}
\end{equation}
for some constant $\lambda$. So using Lemma 4 in \cite{Max05} we can conclude that either $u\equiv 0$, which is impossible as $\|u\|^2_{L^2_\delta(\Omega)}+\|\Tr u\|^2_{L^2(\Sigma)}=1$, or $u>0$ as we wanted to prove.
\end{proof}
\end{proposition}

\begin{lemma}
If $\Omega\subset M$ and $\Sigma\subset \partial M$ are measurable sets, $\lambda_\delta(\Omega,\Sigma)<0$ if and only if $\Yam_b^{q,r}(\Omega,\Sigma)<0$. The result is independent of the indices $\delta$, $q$, $r$ and $b$.
\begin{proof}
In light of Lemma \ref{2l:polynomial}, both quantities are negative if, and only if, there is $u\in W^{1,2}_{\delta^\ast}(\Omega,\Sigma)$ such that $E(u)<0$.
\end{proof}
\end{lemma}
The version of this lemma for the case of positive invariants is a little more involved.
\begin{lemma}
\label{2l:yam_negative_lambda_negative}
If $\Omega\subset M$ and $\Sigma\subset \partial M$ are measurable sets, the following are equivalent:

1. $\Yam_1^{2\bar q, \bar q+1}(\Omega,\Sigma)>0$,

2. $\lambda_\delta(\Omega,\Sigma)>0$ for all $\delta>\delta^\ast$,

3. There is $\delta>\delta^\ast$ such that $\lambda_\delta(\Omega,\Sigma)>0$.
\begin{proof}
First we prove that $1\Rightarrow2$.

If $\delta>\delta^\ast$, by Proposition \ref{2p:minimizer_lambda} implies there is $u\in W^{1,2}_{\delta^\ast}(\Omega,\Sigma)$ such that
\begin{equation}
    \lambda_\delta(\Omega,\Sigma)=\frac{E(u)}{\|u\|^2_{L^2_\delta(\Omega)}+\|\Tr u\|^2_{L^2(\Sigma)}}.
\end{equation}
So if $\lambda_\delta(\Omega,\Sigma)=0$, $E(u)=0$ and, by Lemma \ref{2l:polynomial}, there is $\tilde u\in B^{2\bar q, \bar q+1}_1(\Omega,\Sigma)$ such that $E(\tilde u)=0$, contradicting the assumption. So by Lemma \ref{2l:yam_negative_lambda_negative}, $\lambda_{\delta}(\Omega,\Sigma)>0$.

$2\Rightarrow 3$ is obvious, so now we have to prove $3\Rightarrow 1$.

Assume $\lambda_\delta(\Omega,\Sigma)>0$ and $\Yam^{2\bar q,\bar q+1}_1(\Omega,\Sigma)=0$. Then, there is $\{u_k\}_k\subset B^{2\bar q, \bar q+1}_1(\Omega,\Sigma)$ such that $E(u_k)\to 0$ and, as $\lambda_\delta(\Omega,\Sigma)>0$, $\|u_k\|_{L^{2}_{\delta}(\Omega)},\ \|\gamma u\|_{L^2(\Sigma)}\to 0$. As a consequence, equation (\ref{2eq:gradient_control_by_energy_and_norms}) tells us that $\|\nabla u_k\|_{L^2(\Omega)}\to 0$ and hence the Sobolev inequality (\ref{2eq:sobolev_inequality}) implies $\|u_k\|_{L^{2\bar q}(\Omega)}\to 0$ and, as $u_k\in B^{2\bar , \bar q+1}_1(\Omega,\Sigma)$, $\|\Tr u_k\|_{L^{\bar q+1}(\Sigma)}\to 1$, which contradicts the trace inequality in $\tilde K$. So $\Yam^{2\bar q, \bar q+1}_1(\Omega,\Sigma)>0$ by Lemma \ref{2l:yam_negative_lambda_negative}.
\end{proof}
\end{lemma}

Some properties of these invariants hold as for compact manifolds with boundary, cf. \cite{ST21}.

\begin{proposition}[Monotonicity]
\label{2p:monotonicity}
If $\Omega\subset\tilde\Omega$ and $\Sigma\subset\tilde\Sigma$, then
\begin{equation}
    \Yam^{q,r}_b(\Omega,\Sigma)\geq \Yam^{q,r}_b(\tilde\Omega,\tilde\Sigma)
\end{equation}
and
\begin{equation}
    \lambda_\delta(\Omega,\Sigma)\geq\lambda_\delta(\tilde\Omega,\tilde\Sigma).
\end{equation}
\begin{proof}
The first inequality follows from the fact that $B^{q,r}_b(\Omega,\Sigma)\subset B^{q,r}_b(\tilde\Omega,\tilde\Sigma)$ while the second is a consequence of $W^{1,2}_{\delta^\ast}(\Omega,\Sigma)\subset W^{1,2}_{\delta^\ast}(\tilde\Omega,\tilde\Sigma)$.
\end{proof}
\end{proposition}

\begin{proposition}[Continuity from Above]
\label{2p:continuity_from_above}
If $\{\Omega_k\}_k$ is a decreasing sequence of subsets in $M$ and $\{\Sigma_k\}_k$ is a decreasing sequence of subsets in $\partial M$ such that $\displaystyle\cap_k \Omega_k=\Omega$ and $\displaystyle\cap_k\Sigma_k=\Sigma$, then
\begin{equation}
    \lim_{k\to\infty}\lambda_\delta(\Omega_k,\Sigma_k)=\lambda(\Omega,\Sigma).
\end{equation}
\begin{proof}
By monotonicity, $\{\lambda_\delta(\Omega_k,\Sigma_k)\}_k$ is a non-decreasing sequence and if we call $\Lambda:=\displaystyle\lim_{k\to\infty}\lambda_\delta(\Omega_k,\Sigma_k)$, then $\lambda_\delta(\Omega,\Sigma)\geq\Lambda$.

If $\Lambda=+\infty$, the result is trivial. So assume $\Lambda<+\infty$, hence all the $W^{1,2}_{\delta^\ast}(\Omega_k,\Sigma_k)$ are nontrivial and, given proposition \ref{2p:minimizer_lambda}, for each $k$, there is $u_k\in W^{1,2}_{\delta^\ast}(\Omega_k,\Sigma_k)$ such that:

i. $\|u_k\|^2_{L^2_\delta(\Omega_k)}+\|\Tr u_k\|^2_{L^2(\Sigma_k)}=1$ and

ii. $E(u_k)=\lambda_\delta(u_k)$.

Now, as the sequence is increasing, $\Lambda\geq E(u_k)$ for any $k$, so as in the proof of proposition \ref{2p:minimizer_lambda}, the sequence of the $\{u_k\}_k$ is bounded in $W^{1,2}_{\delta^\ast}(M)$ and there is $u\in W^{1,2}_{\delta^\ast}(M)$ such that:

1. $u_k\rightharpoonup u$ in $W^{1,2}_{\delta^\ast}(M)$,

2. $u_k\to u$ in $L^2_\delta(M)$,

3. $\Tr u_k\to \Tr u$ in $L^2(\partial M)$ and then

4. $\|u\|^2_{L_\delta^2(M)}+\|\Tr u\|^2_{L^2(\partial M)}=1$ and also

5. $E(u)\leq \liminf E(u_k)=\Lambda$.

But 2 and 3 imply that $u_k\to u$ pointwise in $M$ and $\Tr u_k\to \Tr u$ pointwise in $\partial M$, so $u\vert_{(\Omega\cup\Sigma)^c}=0$ and thus $u\in W^{1,2}_{\delta^\ast}(\Omega,\Sigma)$, so, by 5, $\Lambda\geq E(u)\geq \lambda_\delta(\Omega,\Sigma)$, finishing the proof.
\end{proof}
\end{proposition}


\section{The Prescribed Scalar-Mean Curvature Problem}
\label{2sec:prescribed}

In this section we prove our main existence result, establishing the necessary conditions to realize a given pair of functions $R'$ and $H'$ as scalar curvature in $M$ and mean curvature on $\partial M$, respectively, inside a given conformal class. Throughout the section we consider $g$ a $W^{2,p}_{\tau}$ asymptotically euclidean metric on $M$ with $\tau<0$ and $p>\frac{n}{2}$. The assumed regularity of the initial mean curvature on the boundary $H$ and the target one $H'$ is $W^{1-\frac{1}{p}, p}$ as should be induced by $g$ in the whole section, and eventually in our main result, Theorem \ref{2t:general_existence_result}, we will deal with both $R$ and $R'$ in $L^{p}_{\tau-2}(M)$ as well.

In contrast, the intermediate technical lemmata and propositions leading up to the theorem demand stricter conditions on $R$ and $R'$, and we chose at each step to make explicit what is the minimal regularity demanded by each proof to hold. For simplicity, the reader can consider that both $R$ and $R'$ have compact support from Lemma \ref{2l:coercivity_of_F_q_r} to Corollary \ref{2c:uniform_bounds_u_q_r}, an assumption that combined with the $L^p_{\tau-2}$ regularity satisfies the intermediate ones. Lemmata \ref{2l:existence_of_compactly_supported_curvature} and \ref{2l:reduce_curvature} guarantee that such a hypothesis is not too restrictive when we are proving the final theorem.

We will study the functional
\begin{equation}
    F_{q,r}(u)=E(u+1)-\frac{n-2}{2q(n-1)}\int_M R'|u+1|^qdV_g-\frac{n-2}{r}\int_{\partial M}H'|\gamma(u+1)|^rd\sigma_g
\end{equation}
defined over $W^{1,2}_{\delta^\ast}(M,\partial M)$. It will be useful to denote by $Z$ the zero set of $R'$ and by $Z_\partial$ the zero set of $H'$ and define a new norm on $W^{1,2}_{\delta^\ast}(M)$ as
\begin{equation}
    \|u\|^2_{L^2_\delta (M,\partial M)}=\|u\|^2_{L^2_\delta(M)}+\|\gamma u\|^2_{L^2(\partial M)}.
\end{equation}

We start with the technical result below.

\begin{lemma}[Coercivity]
\label{2l:coercivity_of_F_q_r}
Assume $R\in L^1(M)$, $R'\in L^p_{\tau-2}(M)$ and $H'\in W^{1-\frac{1}{p}, p}(\partial M)$. Let $q_0$, $r_0\geq 2$, $q_0<2\bar q$, $r_0<\bar q+1$, $R'\leq 0$, $H'\leq 0$. If $(Z,Z_\partial)$ is a Yamabe positive pair, for all $B\in \mathbb{R}$, there is $K(q_0, r_0, B)>0$ such that if $q$, $r$ are such that $q_0\leq q<2\bar q$ and $r_0\leq r<\bar q+1$ and $u\in W^{1,2}_{\delta^\ast}(M,\partial M)$ satisfies $\|u\|_{L^2_\delta(M,\partial M)}\geq K$, then $F_{q,r}(u)\geq B$.
\begin{proof}
The proof follows the one for Lemma 4.2 in \cite{ST21}, but we didn't manage to shortcut the proof using Section \ref{2sec:gluing}.

For $\epsilon>0$, define the set
\begin{equation}
\begin{array}{ccl}
    A_\epsilon & = & \Big\{u\in W^{1,2}_{\delta^\ast}(M,\partial M), u\geq -1:\int_M|R'|u^2dV_g+\int_{\partial M}|H'|(\Tr u)^2d\sigma_g \\
     & &\leq \epsilon \|u\|^2_{L^2_\delta(M,\partial M)}\left(\int_M |R'|dV_g+\int_{\partial M}|H'|d\sigma_g\right) \Big\}.
\end{array}
\end{equation}
Let $L>0$ be such that $0<L<\lambda_\delta (Z,Z_\partial)$, which exists because $(Z,Z_\partial)$ is Yamabe positive. First, we prove that there is $\epsilon_0<1$ such that if $u\in A_{\epsilon_0}$ then
\begin{equation}
\label{2e:L_condition}
    E(u)\geq L\|u\|^2_{L^2_\delta(M,\partial M)}.
\end{equation}
Assume that it is false. Then if $\epsilon_k$ is a sequence such that $\epsilon_k\to 0$, we can choose $v_k\in A_{\epsilon_k}$, $\|v_k\|_{L^2_\delta(M,\partial M)}=1$ such that
\begin{equation}
    E(v_k)<L.
\end{equation}
Combining equations (\ref{2eq:poincare_inequality}) and (\ref{2eq:gradient_control_by_energy_and_norms}) we get that the sequence $\{v_k\}_{k}$ is bounded in $W^{1,2}_{\delta^\ast}(M)$ so, up to a subsequence, there is $v\in W^{1,2}_{\delta^\ast}(M)$ such that

i. $v_k\rightharpoonup v$ in $W^{1,2}_{\delta^\ast}(M)$;

ii. $v_k\to v$ in $L^2_\delta(M)$ by the compactness of the embedding $W^{1,2}_\delta(M)\to L^2_\delta(M)$;

iii. $\Tr v_k \to \Tr v$ in $L^2(\partial M)$ by the compactness of the trace map.

From ii and iii, $\|v\|_{L^2_\delta (M,\partial M)}=1$. Also, as proven in \cite{ST21}, the map $u\to \int_{\partial M}H'u^2d\sigma_g$ is continuous in $L^2(\partial M)$, hence, by iii,
\begin{equation}
    \int_{\partial M}H'(\Tr v_k)^2d\sigma_g\to \int_{\partial M}H'(\Tr v)^2d\sigma_g.
\end{equation}
Moreover, by Lemma 3.1 in \cite{DM15},
\begin{equation}
    \int_{M}R' v_k^2d\sigma_g\to \int_{M}R'v^2d\sigma_g.
\end{equation}
Now, by choice of the $v_k$'s
\begin{equation}
    0\leq \int_M |R'|v_k^2dV_g+\int_{\partial M}|H'|(\gamma v_k)^2d\sigma_g\leq \epsilon_k ||v_k||^2_{L^2(M,\partial M)}\left(\int_M|R'|dV_g+\int_{\partial M}|H'|d\sigma_g\right)\to 0,
\end{equation}
hence
\begin{equation}
    \int_M R' v^2dV_g+\int_{\partial M}H'(\Tr v)^2d\sigma_g=0
\end{equation}
and $v\in W^{1,2}_{\delta^\ast}(Z.Z_\partial)$. Besides, by Corollary 3.2 of \cite{DM15}
\begin{equation}
    E(v)\leq \lim E(v_k)\leq L<\lambda_\delta(Z,Z_\partial),
\end{equation}
contradiction.

With this first result proven, we have two possibilities next.

If $u\in A_{\epsilon_0}$, as $R',H'\leq 0$:
\begin{equation*}
\begin{array}{ccl}
    F_{q,r}(u) & \geq & E(u+1) \\
     &= & E(u) + \int_M R((u+1)^2-u^2)dV_g +\int_{\partial M}H((\gamma(u+1))^2-(\gamma u)^2)d\sigma_g\\
     &=& E(u) + \int_M R(2u+1)dV_g +\int_{\partial M}H(2\gamma u+1)d\sigma_g\\
     &\geq& E(u)-\int_M|R|\left(1+\epsilon u^2+\frac{1}{\epsilon}\right)dV_g-\int_{\partial M}|H|\left(1+\epsilon(\Tr u)^2+\frac{1}{\epsilon}\right)d\sigma_g,
\end{array}
\end{equation*}
because $\left(\sqrt{\epsilon}v-\frac{1}{\sqrt{\epsilon}}\right)^2=\epsilon v^2+\frac{1}{\epsilon}-2v\geq 0$ for any $\epsilon>0$.
\begin{equation*}
    \begin{array}{ccl}
       F_{q,r}(u) & \geq & (1-\epsilon)E(u)+\epsilon\left(\int_{M}|\nabla u|^2dV_g-\frac{5n-6}{4(n-1)}\int_M|R|u^2dV_g-\frac{n}{2}\int_{\partial M}|H|(\Tr u)^2d\sigma_g\right) \\
         & & -\left(1+\frac{1}{\epsilon}\right)\int_M|R|dV_g-\left(1+\frac{1}{\epsilon}\right)\int_{\partial M}|H|d\sigma_g.
    \end{array}
\end{equation*}
Now, Lemma \ref{2l:bounds_on_curvature_terms} does not specifically depend on $R$ and $H$ being the curvatures, but only on their regularity. Hence, the inequalities hold for $|R|$ and $|H|$ as well and we have that there is, in the notation of Lemma \ref{2l:bounds_on_curvature_terms}, $C>\max\left\{K_{\frac{4(n-1)}{3(5n-6)}}, K_{\frac{2}{3n}}\right\}$ such that
\begin{equation*}
    \begin{array}{ccl}
        F_{q,r}(u) &\geq &(1-\epsilon)E(u)+\frac{\epsilon}{3}\int_M|\nabla u|^2dV_g-\epsilon C\|u\|^2_{L^2_\delta(M)}-\epsilon C\|\Tr u\|^2_{L^2(\partial M)}\\
        &&-\left(1+\frac{1}{\epsilon}\right)\left(\int_M |R| dV_g+\int_{\partial M}|H|d\sigma_g\right)  \\
        &\geq & (L-\epsilon(L+C))\|u\|^2_{L^2_\delta(M,\partial M)}+\frac{\epsilon}{3}\int_M|\nabla u|^2dV_g-\left(1+\frac{1}{\epsilon}\right)\left(\int_M |R| dV_g+\int_{\partial M}|H|d\sigma_g\right),
    \end{array}
\end{equation*}
because $u\in A_{\epsilon_0}$. The result follows by choosing $\epsilon<\frac{L}{L+C}$ because $R\in L^1(M)$ and $H\in L^1(\partial M)$.

On the other hand, assume $u\not\in A_{\epsilon_0}$. Thus
\begin{equation}
    \int_M|R'|u^2dV_g+\int_{\partial M}|H'|(\Tr u)^2d\sigma_g\geq \epsilon_0\|u\|^2_{L^2_\delta(M,\partial M)}\left(\int_M|R'|dV_g+\int_{\partial M}|H'|d\sigma_g\right).
\end{equation}
In this case, since $R'\leq 0$, $H'\leq 0$
\begin{equation*}
\begin{array}{ccl}
    F_{q,r}(u)&=&E(u+1)+\frac{n-2}{2q(n-1)}\int_M|R'||u+1|^qdV_g+\frac{n-2}{r}\int_{\partial M}|H'||\gamma u+1|^rd\sigma_g\\
    &\geq & E(u+1)-\frac{n-2}{2q(n-1)}\int_M|R'|dV_g-\frac{n-2}{r}\int_{\partial M}|H'|d\sigma_g\\
    &&+\frac{n-2}{2q(n-1)}\int_M|R'||u|^qdV_g+\frac{n-2}{r}\int_{\partial M}|H'||\gamma u|^rd\sigma_g,
\end{array}
\end{equation*}
because $a\geq -1$ implies $(a+1)^p\geq |a|^p-1$ for $p\geq 1$.

Using H\"older's inequality as in the proof of Lemma 4.2 in \cite{ST21} one gets that:
\begin{equation}
\begin{array}{ccl}
    F_{q,r}(u) &\geq& E(u+1)-\frac{n-2}{2q(n-1)}\int_M|R'|dV_g-\frac{n-2}{r}\int_{\partial M}|H'|d\sigma_g \\
     & &+A\left(\int_M|R'|u^2dV_g\right)^{\frac{q}{2}}+B\left(\int_{\partial M}|H'|(\Tr u)^2d\sigma_g\right)^{\frac{r}{2}}
\end{array}
\end{equation}
with $A=\displaystyle\min_{q\in [q_0,2\bar q]}\frac{n-2}{2q(n-1)}\left(\int_M |R'|dV_g\right)^{1-\frac{q}{2}}$ and $B=\displaystyle\min_{r\in [r_0,\bar q+1]}\frac{n-2}{r}\left(\int_{\partial M}|H'|d\sigma_g\right)^{1-\frac{r}{2}}$.

On the other hand,
\begin{equation*}
    \begin{array}{ccl}
        E(u+1) & = & \int_M|\nabla u|^2dV_g+\frac{n-2}{4(n-1)}\int_M R(u+1)^2dV_g+\frac{n-2}{2}\int_{\partial M}H(\Tr u+1)^2d\sigma_g \\
        & \geq & \int_M|\nabla u|^2dV_g - \frac{n-2}{2(n-1)}\int_M|R|(u^2+1)dV_g-(n-2)\int_{\partial M}|H|\left((\Tr u)^2+1\right)d\sigma_g\\
        & \geq & \frac{1}{3}\int_M|\nabla u|^2dV_g - \frac{n-2}{2(n-1)}\int_M|R|dV_g-(n-2)\int_{\partial M}|H|d\sigma_g-C\|u\|^2_{L^2_\delta(M,\partial M)}
    \end{array}
\end{equation*}
as in the previous case. The result then follows from similar calculations as in Case 1 in the proof of Lemma 4.2 in \cite{ST21}, as $R\in L^1(M)$ and $H\in W^{1-\frac{1}{p}, p}(\partial M)\hookrightarrow L^1(\partial M)$.

\end{proof}
\end{lemma}

Among other things, the coercivity result above guarantees that if $W^{1,2}_{\delta^\ast}(M,\partial M)$ is nonempty, the functional $F_{q,r}$ is bounded from below, and we can look for minimizers.

\begin{proposition}
\label{2p:existence_subcritical_minimizers}
If $R$, $R'$, $H'$, $q_0$ and $r_0$ are as in the previous lemma, $H'$ bounded, $2\bar q>q\geq q_0$, $\bar q+1>r\geq r_0$, $q>r$ and $(Z,Z_\partial)$ is a Yamabe positive pair, then there is $u_{q,r}>-1$ in $W^{1,2}_{\delta^\ast}(M)$ that minimizes $F_{q,r}$.

If in addition $R$ and $R'$ have compact support, $u_{q,r}\in W^{2,p}_\delta(M)$ for any $\delta\in(2-n,0)$.
\begin{proof}
Let $\{u_k\}_k\subset W^{1,2}_{\delta^\ast}(M)$. Notice that, for any $u$, $F(u)=F(-u-2)$, so we can choose $u_k\geq -1$ for all $k$.

Since the sequence is minimizing, we can assume up to a subsequence, that
\begin{equation}
\label{2e:energy_bounded_by_zero_energy}
    F_{q,r}(u_k)\leq F_{q,r}(0):=B,
\end{equation}
because $0\in W^{1,2}_{\delta^\ast}(M,\partial M)$. Then the coercivity implies that there is $K$ not depending on $q$ and $r$ - if we assume our sequences satisfy condition (\ref{2e:energy_bounded_by_zero_energy}) - such that
\begin{equation}
    \|u_k\|^2_{L^2_\delta(M,\partial M)}\leq K.
\end{equation}
Also
\begin{equation}
    B\geq F_{q,r}(u_k)\geq E(u_k+1)
\end{equation}
and since $\nabla u_k=\nabla(u_k+1)$, we finally get that $\|u_k\|_{W^{1,2}_{\delta^\ast}(M)}$ is uniformly bounded, and there is $u_{q,r}\in W^{1,2}_{\delta^\ast}(M)$ a weak limit for the sequence. As in the proof of the previous lemma we have that
\begin{equation}
    \int_M R (u_k+1)^2 dV_g\to \int_M R(u_{q,r}+1)^2 dV_g
\end{equation}
and
\begin{equation}
    \int_{\partial M}H(\Tr u_k+1)^2d\sigma_g\to \int_{\partial M} H(\Tr u_{q,r}+1)^2d\sigma_g,
\end{equation}
hence, $E(u_{q,r}+1)\leq \liminf E(u_k+1)$.

Similarly, $\Tr u_k\to \Tr u_{q,r}$ in $L^r(\partial M)$ and $H'$ is bounded, so
\begin{equation}
    \int_{\partial M}H'(\Tr (u_k+1))^rd\sigma_g\to \int_{\partial M}H'(\Tr (u_{q,r}+1))^rd\sigma_g.
    \end{equation}
Moreover, as $W^{1,2}_{\delta^\ast}(\Omega)\hookrightarrow L^q_\delta(\Omega)$, $u_k\to u$ in $L^q_\delta(\Omega)$ for any $\delta>\delta^\ast$ and, as in Lemma 4.6 in \cite{DM15},
    \begin{equation}
        \int_M R' (u_k+1)^q dV_g\to \int_M R'(u_{q,r}+1)^q dV_g
    \end{equation}
We have that $F_{q,r}(u)\leq F_{q,r}(u_k)$ for all $k$, so $u$ is our minimizer. Furthermore, by our choice of the $u_k$'s, as $u_k\to u$ pointwise, $u\geq -1$.

However if $u$ is a minimizer, it is a weak solution to the equation

\begin{equation}
    \label{2eq:subcritical_problem}
\begin{cases}
    -\Delta u+\frac{n-2}{4(n-1)}R(u+1)=\frac{n-2}{4(n-1)}R'(u+1)^{q-1},\ in\ M,\\
    \Tr \partial_\nu u+\frac{n-2}{2}H(\Tr u+1)=\frac{n-2}{2}H'(\Tr u+1)^{r-1},\ on\ \partial M,
\end{cases}
\end{equation}
hence Lemma 4 from \cite{Max05} implies that either $u+1\equiv 0$ or $u>-1$. But $u\to 0$ towards infinity, as $u\in W^{1,2}_{\delta^\ast}(M)$, so $u>-1$.

Now, assume that $R$ and $R'$ have compact support, that $K$ and $\tilde K$ are compact sets as in section \ref{2sec:gluing} and that $R$ and $R'$ vanish outside $K$.

Using Proposition 1 in \cite{Max05}, we can get the following estimate when considering the operator $\mathcal{P}=(-\Delta, \Tr\partial_{\nu})$:
\begin{equation}
    \|u_{q,r}\|_{W^{2,p}_\delta(M)}\leq C\left(\|\Delta u_{q,r}\|_{L^p_{\delta-2}(M)}+\|\Tr \partial_\nu u_{q,r} \|_{W^{1-\frac{1}{p}, p}(\partial M)}\right)
\end{equation}
for some $C>0$. But $u_{q,r}$ satisfies equation (\ref{2eq:subcritical_problem}). Using the triangle inequality we get
\begin{equation*}
\begin{array}{ccl}
    \|u_{q,r}\|_{W^{2,p}_\delta(M)} & \leq &  \frac{C(n-2)}{4(n-1)}\left(\|R(u_{q,r}+1)\|_{L^p_{\delta-2}(M)}+\|R'(u_{q,r}+1)^{q-1}\|_{L^p_{\delta-2}(M)}\right)\\
    &&+\frac{C(n-2)}{2}\left(\|H\Tr (u_{q,r}+1) \|_{W^{1-\frac{1}{p}, p}(\partial M)}+\|H'\Tr(u_{q,r}+1)^{r-1}\|_{W^{1-\frac{1}{p}, p}(\partial M)}\right)\\
   & = &  \frac{C(n-2)}{4(n-1)}\left(\|R(u_{q,r}+1)\|_{L^p_{\delta-2}(\tilde K)}+\|R'(u_{q,r}+1)^{q-1}\|_{L^p_{\delta-2}(\tilde K)}\right)\\
    &&+\frac{C(n-2)}{2}\left(\|H\Tr (u_{q,r}+1) \|_{W^{1-\frac{1}{p}, p}(\partial M)}+\|H'\Tr(u_{q,r}+1)^{r-1}\|_{W^{1-\frac{1}{p}, p}(\partial M)}\right),
\end{array}
\end{equation*}
because the supports of $R$ and $R'$ are contained in $\tilde K$. But $\tilde K$ is compact, hence the weighted and non-weighted norms are equivalent in $\tilde K$ and we can drop the $\delta-2$ in the indices. Also, using the fact that $(a+b)^s\leq 2^{s-1}(a^s+b^s)$ for $s\geq 1$ to split the non-linear terms and the triangle inequality, we can find $C_2>0$ such that
\begin{equation*}
\begin{array}{ccl}
    \|u_{q,r}\|_{W^{2,p}_\delta(M)} & \leq & C_2\left(\|Ru_{q,r}\|_{L^p(\tilde K)}+\|R\|_{L^p(\tilde K)}+\|R'u_{q,r}^{q-1}\|_{L^p(\tilde K)}+\|R'\|_{L^p(\tilde K)}\right)\\
    &&+C_2\left(\|H\Tr u_{q,r} \|_{W^{1-\frac{1}{p}, p}(\partial M)}+\|H\|_{W^{1-\frac{1}{p}, p}(\partial M)}+\|H'\Tr(u_{q,r}+1)^{r-1}\|_{W^{1-\frac{1}{p}, p}(\partial M)}\right).
\end{array}
\end{equation*}
Finally, as the function spaces of the norms on the right hand side of the inequality are with respect to a compact manifold with boundary, we can use the techniques described in section \ref{2sec:gluing} and the bootstrap procedure described in the proof of Lemma 4.5 of \cite{ST21} to get the regularity needed. Notice only that to control the last term one eventually needs to control the term $\|u_{q,r}+1\|_{W^{1,2}(\tilde K)}$, but
\begin{equation}
    \|u_{q,r}+1\|^2_{W^{1,2}(\tilde K)}=\|\nabla u_{q,r}\|^2_{L^2(\tilde K)}+\|u_{q,r}+1\|^2_{L^2(\tilde K)}\leq \|\nabla u_{q,r}\|^2_{L^2(\tilde K)}+2\|u_{q,r}\|^2_{L^2(\tilde K)}+2V_g(\tilde K),
\end{equation}
and that is easy to handle as $\tilde K$ has finite measure and $\|u_{q,r}\|_{W^{1,2}_{\delta^\ast}(M)}$ is uniformly bounded as mentioned in the commentary following this proof and equation (\ref{2eq:uniform_1,2_bound_for_the_u_qr}) below. 
\end{proof}
\end{proposition}

Notice that, as $B$ and $K$ in the previous proof do not depend on $q$ and $r$, the uniform bounds on each sequence $\{\|u_k\|_{W^{1,2}_{\delta^\ast}(M)}\}_k$ does not depend on $q$ and $r$ either. As a consequence, there is $C>0$ that does not depend on $q$ or $r$ such that
\begin{equation}
\label{2eq:uniform_1,2_bound_for_the_u_qr}
    \|u_{q,r}\|^2_{W^{1,2}_{\delta^\ast}(M)}=\|\nabla u_{q,r}\|^2_{L^2(M)}+\|u_{q,r}\|^2_{L^2_{\delta^\ast}(M)}<C.
\end{equation}

The same bootstrap procedure just used shows a higher regularity for $u_{q,r}$ when we guarantee an {\it{a priori}} bound for $\|u_{q,r}\|_{L^Q(\tilde K)}$ with $Q>2\bar q$, as the one we will prove the following lemma.

\begin{lemma}
\label{2l:bound_for_Q_norm}
For any compact set $K\subset M\cup\partial M$, there is $C_{K}$ a uniform bound such that $\|u_{q,r}\|_{L^Q(K)}<C_{K}$ for any $q$, $r$ subcritical, for some $Q>2\bar q$, $Q$ depending on $K$.

\begin{proof}
This proof adapts the ideas in the proof of Lemma 4.8 in \cite{DM15} to manifolds with boundary.

If $u_{q,r}$ is a minimizer of $F_{q,r}$, we define the following auxiliary functions for a given $\delta>0$:
\begin{equation}
    \begin{cases}
        w=(1+u_{q,r})^{1+\delta}\\
        \nu=\xi_{\tilde K}^2(1+u_{q,r})^{1+2\delta}
    \end{cases}
\end{equation}
with $\xi_{\tilde K}$ a smooth non-negative function with compact support that is 1 in $K$ constructed as in section \ref{2sec:gluing}. Notice that $(u_{q,r}+1)\nu=\xi_{\tilde K}^2w^2$. Our first goal is to find an uniform bound to
\begin{equation}
    \|\nabla(\xi_{\tilde K}w)\|^2_{L^2(M)}=\|\xi_{\tilde K}\nabla w+w\nabla \xi_{\tilde K}\|^2_{L^2(M)}.
\end{equation}
Using the triangle inequality and the fact that for positive numbers $a$, $b$, $2ab\leq 2a^2+2b^2$ we can work that estimate up to
\begin{equation}
    \|\nabla(\xi_{\tilde K}w)\|^2_{L^2(M)}\leq 3\left(\|\xi_{\tilde K}\nabla w\|^2_{L^2(M)}+\|w\nabla\xi_{\tilde K}\|^2_{L^2(M)}\right).
\end{equation}

The second term is easy to deal with. As $\xi_{\tilde K}$ is supported in $\tilde K$ and is smooth:
\begin{equation*}
\begin{array}{ccl}
    \|w\nabla\xi_{\tilde K}\|^2_{L^2(M)} & = & \|w\nabla\xi_{\tilde K}\|^2_{L^2(\tilde K)}\\
    &\leq & \left(\max_{\tilde K}|\nabla\xi_{\tilde K}|^2\right)\|(1+u_{q,r})^{1+\delta}\|^2_{L^2(\tilde K)}\\
    &=&\left(\max_{\tilde K}|\nabla\xi_{\tilde K}|^2\right)\|1+u_{q,r}\|^{2(1+\delta)}_{L^{2(1+\delta)}(\tilde K)}
\end{array}
\end{equation*}
and because $\tilde K$ is compact and using the triangle inequality
\begin{equation}
    \|w\nabla\xi_{\tilde K}\|^2_{L^2(M)}\leq 2^{2(1+\delta)} \left(\max_{\tilde K}|\nabla\xi_{\tilde K}|^2\right)\left(V_g(\tilde K)+\|u_{q,r}\|^{2(1+\delta)}_{L^{2(1+\delta)}(\tilde K)}\right).
\end{equation}
So if we force as a first restriction that $\delta$ is so small that $2(1+\sigma)<2\bar q$, because $\tilde K$ is compact we can use inequality (\ref{2eq:uniform_1,2_bound_for_the_u_qr}), the usual Sobolev embedding results in compact manifolds and the fact that
\begin{equation}
\|u_{q,r}\|^2_{W^{1,2}(\tilde K)}\leq \max_{\tilde K}\left(1+\rho^{-2\delta^\ast-n}\right)\|u_{q,r}\|^2_{W^{1,2}_{\delta^\ast}(\tilde K)} \leq \max_{\tilde K}\left(1+\rho^{-2\delta^\ast-n}\right)\|u_{q,r}\|^2_{W^{1,2}_{\delta^\ast}(M)}
\end{equation}
to get an uniform bound on $\|w\nabla\xi_{\tilde K}\|^2_{L^2(M)}$ that does not depend on $q$ or $r$.

With respect to the first term,
\begin{equation}
    \langle \nabla u_{q,r}, \nabla \nu\rangle=\frac{1+2\delta}{(1+\delta)^2}\xi^2_{\tilde K}\langle\nabla w,\nabla w\rangle+(1+u_{q,r})^{2\delta+1}\langle \nabla u_{q,r},\nabla(\xi_{\tilde K})^2\rangle,
\end{equation}
so
\begin{equation}
    \langle \xi_{\tilde K}\nabla w,\xi_{\tilde K}\nabla w\rangle=\frac{(1+\delta)^2}{1+2\delta}\left(\langle \nabla u_{q,r}, \nabla \nu\rangle-(1+u_{q,r})^{2\delta+1}\langle\nabla u_{q,r}, \nabla(\xi_{\tilde K})^2\rangle\right).
\end{equation}

Again, it is easy to control the last term. First, by the Cauchy-Schwartz inequality
\begin{equation}
    (1+u_{q,r})^{2\delta+1}\langle\nabla u_{q,r}, \nabla(\xi_{\tilde K})^2\rangle\leq (1+u_{q,r})^{2\delta+1}|\nabla u_{q,r}||\nabla(\xi_{\tilde K})^2|
\end{equation}
and as a result, by the same inequality applied to integrals, and using the fact that $\xi_{\tilde K}$ is smooth and $\tilde K$ is compact
\begin{equation}
    \int_{M}(1+u_{q,r})^{2\delta+1}\langle\nabla u_{q,r}, \nabla(\xi_{\tilde K})^2\rangle dV_g\leq \|(1+u_{q,r})\|^{2\delta+1}_{L^{4\delta+2}(\tilde K)}\max_{\tilde K}|\nabla(\xi_{\tilde K})^2|\|\nabla u_{q,r}\|_{L^2(\tilde K)}
\end{equation}
and inequality (\ref{2eq:uniform_1,2_bound_for_the_u_qr}) guarantees an uniform bound for $\|\nabla u_{q,r}\|_{L^2(\tilde K)}$ and for $\|u_{q,r}\|^{4\delta+2}_{L^{4\delta+2}(\tilde K)}$ as long as we force as a second restriction on $\delta$ that $4\delta+2<2\bar q$. As a result, as $\tilde K$ is compact, $\|1+u_{q,r}\|_{L^{4\delta+2}(\tilde K)}$ is also bounded uniformly with respect to $q$ and $r$ and hence the term we are addressing here, $\int_{M}(1+u_{q,r})^{2\delta+1}\langle\nabla u_{q,r}, \nabla(\xi_{\tilde K})^2\rangle dV_g$, is as well.

Finally, as $u_{q,r}$ is a minimizer of $F_{q,r}$, it is a weak solution of equation (\ref{2eq:subcritical_problem}). So testing $u_{q,r}$ against $\nu$ we have that
\begin{equation*}
\begin{array}{ccl}
\int_M\langle \nabla u_{q,r}, \nabla \nu \rangle dV_g&=&\frac{n-2}{4(n-1)}\int_M \left(R'(u_{q,r}+1)^{q-2}-R\right)(u_{q,r}+1)\nu dV_g\\
&&+\frac{n-2}{2}\int_{\partial M} \left(H'(\Tr u_{q,r})^{r-2}-H \right)\Tr(( u_{q,r}+1)(\nu)) d\sigma_g\\
&=&\frac{n-2}{4(n-1)}\int_M \left(R'(u_{q,r}+1)^{q-2}-R\right)\xi_{\tilde K}^2w^2 dV_g\\
&&+\frac{n-2}{2}\int_{\partial M} \left(H'(\Tr u_{q,r})^{r-2}-H \right)\Tr(\xi_{\tilde K}^2w^2) d\sigma_g\\
    &\leq & -\frac{n-2}{4(n-1)}\int_M R\xi_{\tilde K}^2w^2 dV_g-\frac{n-2}{2}\int_{\partial M} H \Tr(\xi_{\tilde K}^2w^2) d\sigma_g\\
    &\leq&  \frac{n-2}{4(n-1)}\left\vert\int_M R\xi_{\tilde K}^2w^2 dV_g\right\vert+\frac{n-2}{2}\left\vert\int_{\partial M} H \Tr(\xi_{\tilde K}^2w^2) d\sigma_g\right\vert,
\end{array}
\end{equation*}
the first inequality because $R'\leq 0$ and $H'\leq 0$.

Using the estimates obtained in Lemma 4.1 we can proceed to get that for any $\epsilon>0$, there is $C_\epsilon$ such that for some $\sigma$
\begin{equation}
\begin{array}{ccl}
    \int_{M}\langle\nabla u_{q,r}, \nabla \nu\rangle dV_g & \leq & \epsilon\left( \|\nabla(\xi_{\tilde K}w)\|^2_{L^2(M)}+\|\gamma(\xi_{\tilde K}w)\|^2_{W^{\frac{1}{2}, 2}(\partial M)}\right)\\
     & &+ C_\epsilon \left(\|\xi_{\tilde K}w\|^2_{L^2_\sigma(M)}+\|\gamma(\xi_{\tilde K}w)\|^2_{L^{2}(\partial M)}\right).
\end{array}
\end{equation}
Notice the term multiplying $C_\epsilon$ can be bounded uniformly with respect to $q$ and $r$ with the same choice of $\delta$ as before, since $\|\xi_{\tilde K}w\|_{L^2_\sigma(M)}$ is equivalent to $\|w\|_{L^2(\tilde K)}$ and the trace inequality allows us to control $\|\Tr(\xi_{\tilde K}w)\|^2_{L^2(\partial M)}$, which is equivalent to $\|\Tr u_{q,r}\|^{2\delta+2}_{L^{2\delta+2}(\partial M)}$, by $\|u_{q,r}\|_{W^{1,2}(M)}$. Summing everything up, we have that there is $D_1>0$, independent of $q$ and $r$, such that
\begin{equation}
    \|\nabla(\xi_{\tilde K}w)\|^2_{L^2(M)}\leq 3\frac{(1+\delta)^2}{1+2\delta}\epsilon\left( \|\nabla(\xi_{\tilde K}w)\|^2_{L^2(M)}+\|\gamma(\xi_{\tilde K}w)\|^2_{W^{\frac{1}{2}, 2}(\partial M)}\right)+D_1.
\end{equation}
Now, using the trace inequality followed by the Sobolev inequality in $\tilde K$, we know there is $D_2>0$ such that
\begin{equation}
    \|\nabla(\xi_{\tilde K}w)\|^2_{L^2(M)}\leq 3\frac{(1+\delta)^2}{1+2\delta}(1+D_2)\epsilon \|\nabla(\xi_{\tilde K}w)\|^2_{L^2(M)}+D_1.
\end{equation}
and by choosing $\epsilon$ small enough, we can rearrange it into
\begin{equation}
    \|\nabla(\xi_{\tilde K}w)\|^2_{L^2(M)}<D
\end{equation}
for some $D>0$ that does not depend on $q$ or $r$.

Hence, from the Sobolev inequality, $\|w\|_{L^{2\bar q}(\tilde K)}=\|u_{q,r}^{1+\delta}\|_{L^{2\bar q}(\tilde K)}$ is bounded uniformly with respect to $q$ and $r$ and the result follows for $Q=2\bar q(1+\delta)$.
\end{proof}
\end{lemma}

\begin{corollary}
\label{2c:uniform_bounds_u_q_r}
If $R$, $R'$ have compact support, under the hypothesis and the naming of Proposition \ref{2p:existence_subcritical_minimizers}, given $P>\frac{n}{2}$, $\delta\in(2-n, 0)$, $u_{q,r}\in W^{2,P}_\delta(M)$, $\Tr u_{q,r}\in W^{2,R}(\partial M)$ for some $R>\frac{n-1}{2}$ and there are uniform bounds $C$, $K$ independent of $q\in(2,2\bar q)$ and $r\in (2,\bar q+1)$ such that
\begin{equation}
\begin{cases}
\|u_{q,r}\|_{W^{2,P}_\delta(M)}\leq C,\\
\|\Tr u_{q,r}\|_{W^{2,R}(\partial M)}\leq K.
\end{cases}    
\end{equation}
\begin{proof}
If $\tilde K$ is as in the proof of Proposition \ref{2p:existence_subcritical_minimizers}, the bootstrap mechanism described in the same proposition applied a finite number of times guarantees that there are constants $C_1$, $C_2$ such that
\begin{equation}
    \|u_{q,r}\|_{W^{2,P}_\delta(M)}\leq C_1\|u_{q,r}\|_{L^Q(\tilde K)}+C_2,
\end{equation}
and Lemma \ref{2l:bound_for_Q_norm} guarantees the uniform bound on $\|u_{q,r}\|_{W^{2,P}_\delta(M)}$.

The bound on $\|\Tr u_{q,r}\|_{W^{2,R(\partial M)}}$ then follows from the trace inequality applied to $\xi_{\tilde K}u_{q,r}$.
\end{proof}
\end{corollary}

Up until now we have developed tools that provide conditions for a conformal equivalence between metrics with scalar curvatures $R$ and $R'$ in $M$ as long as both have compact support. To deal with the general case we need a couple of intermediate results to navigate between metrics which corresponding curvature do not have compact support and those that have.

We start with one direction: we always have a metric inside a given conformal class that has a scalar curvature $\tilde R$ of compact support inside $M$. The proof is the same as for the case of asymptotically euclidean manifolds without boundary presented for Lemma 4.4 in \cite{DM15}, as the proof relies only with spaces of functions restricted to the ends.

\begin{lemma}
\label{2l:existence_of_compactly_supported_curvature}
Assume $(M\cup\partial M, g)$ is an asymptotically Euclidean manifold with boundary and metric $g\in W^{2,p}_\tau(M)$. Then there is a metric $\tilde g$ conformal to $g$ such that if $\tilde R$ is the scalar curvature on $M$ induced by $\tilde g$, $\tilde R\equiv 0$ outside of a compact set.
\end{lemma}

We will eventually reduce the problem to finding a suitable metric with a given curvature $\tilde R'$ of compact support that can be realized starting from a given metric $R$. Our goal is to find $\tilde R'\geq R'$ due to the following result saying we can decrease the curvature through a conformal transformation.

\begin{lemma}
\label{2l:reduce_curvature}
Assume $(M\cup\partial M, g)$ is an asymptotically Euclidean manifold with boundary, $g\in W^{2,p}_\tau(M)$. If $\tilde R\leq R$ and $\tilde H\leq H$, $\tilde R\in L^{ p}_{\tau-2}(M)$, $\tilde H\in W^{k-1-\frac{1}{p},p}(\partial M)$, then there is $\tilde g$ a metric in the same conformal class of $g$ that induces $\tilde R$ as the scalar curvature on $M$ and $\tilde H$ as the mean curvature on $\partial M$.
\begin{proof}

We want solutions to the equation:
\begin{equation}
\label{2eq:equation_to_reduce_curvature}
    \begin{cases}
    -\Delta u+\frac{n-2}{4(n-1)}R (u+1)=\frac{n-2}{4(n-1)}\tilde R (u+1)^{2\bar q-1},\ in\ M,\\
\Tr\partial_\nu u+\frac{n-2}{2}H\Tr (u+1)=\frac{n-2}{2}\tilde H(\Tr (u+1))^{\bar q},\ on\ \partial M,
\end{cases}
\end{equation}
and try to find a solution $u>-1$. If there is such $u$, $(u+1)^{2\bar q-2}$ is the desired conformal factor.

Notice $\tilde R\leq R$ and $\tilde H\leq H$, $u_+\equiv 0$ is a supersolution in the terms of section 3 of \cite{Max05}. To be able to use a barrier method, let us try to find a subsolution as well.

Define the functions $\tilde R'=\min\{0, \tilde R\}\leq \tilde R$ in $M$ and $\tilde H'=\min\{0, \tilde H\}\leq \tilde H$ and consider the new equation
\begin{equation}
    \begin{cases}
    -\Delta v+\frac{n-2}{4(n-1)}(R-\tilde R') v=-\frac{n-2}{4(n-1)}(R-\tilde R'),\ in\ M,\\
\Tr\partial_\nu v+\frac{n-2}{2}(H-\tilde H')\Tr v=-\frac{n-2}{2}\tilde (H-\tilde H'),\ on\ \partial M.
\end{cases}
\end{equation}

As $R\geq \tilde R\geq \tilde R'$, and likewise $H\geq \tilde H'$, the coefficients of the linear terms on both equations are non-negative, hence the equation has a solution $v\in W^{k,p}_{\tau}(M)$ by Proposition 1 in \cite{Max05} and by Lemma 2 in the same paper, $v\leq 0$. On the other hand, we have that $1+v$ satisfies
\begin{equation}
    \begin{cases}
    -\Delta (1+v)+\frac{n-2}{4(n-1)}(R-\tilde R')(1+v)=0,\ in\ M,\\
\Tr\partial_\nu (1+v)+\frac{n-2}{2}(H-\tilde H')\Tr (1+v)=0,\ on\ \partial M,
\end{cases}
\end{equation}
so by Lemma 4 in the same same paper, $1+v>0$ everywhere, as $1+v\to 1$ at infinity. Summing up $0<1+v\leq 1$.

Finally, comparing with equation (\ref{2eq:equation_to_reduce_curvature}), we have that, for $v$:
\begin{equation}
\begin{cases}
-\Delta v+\frac{n-2}{4(n-1)}R(v+1)=\frac{n-2}{4(n-1)}\tilde R'(v+1)\leq \frac{n-2}{4(n-1)}\tilde R'(v+1)^{2\bar q-1}\leq \frac{n-2}{4(n-1)} \tilde R(v+1)^{2\bar q-1},\\
\Tr\partial_{\nu}v+\frac{n-2}{2}H\Tr(v+1)=\frac{n-2}{2}\tilde H'\Tr (v+1)\leq \frac{n-2}{2}\tilde H'(\Tr (v+1))^{\bar q-1}\leq \frac{n-2}{2}\tilde H(\Tr (v+1))^{\bar q-1},
\end{cases}
\end{equation}
the first inequality on each line because $1+v\leq 1$, $\tilde R'$ and $\tilde H'$ are non-positive and $\bar q>1$, the second because $\tilde R\geq \tilde R'$ and $\tilde H\geq \tilde H'$ and $1+v>0$. So $v$ is our subsolution.

Finally, Proposition 2 in \cite{Max05} guarantees the existence of a solution $u\in W^{k,p}_{\tau}(M)$ of (\ref{2eq:equation_to_reduce_curvature}), satisfying $-1<v\leq u\leq 0$ and as $1+u\to 1$ at infinity, and the metric $\tilde g=(1+u)^{2\bar q-2}g$ is the one we are looking for.
\end{proof}
\end{lemma}

We can finally prove the existence result.

\begin{theorem}
\label{2t:general_existence_result}
Let $(M,\partial M, g)$ be an asymptotically euclidean Riemannian manifold with boundary, $g\in W^{2,p}_\tau(M)$, $p>\frac{n}{2}$, $\tau<0$. If $R'\leq 0$, $R'\in L^{ p}_{\tau-2}(M)$ and $H'\leq 0$, $H'\in W^{1-\frac{1}{p}, p}(\partial M)$, then there is $g'$ a metric in the same conformal class as $g$ realizing $R'$ as the scalar curvature in $M$ and $H'$ as the mean curvature on $\partial M$ associated to $g'$ if, and only if, $\Yam(Z,Z_\partial)>0$.
\begin{proof}
First, we prove the existence of $g'$ assuming $\Yam(Z,Z_\partial)>0$.

Notice that, by Lemma \ref{2l:existence_of_compactly_supported_curvature}, there is a metric $\tilde g$ conformal to $g$ that realizes $\tilde R$, a scalar curvature in $M$ with compact support. Call $\tilde H$ the mean curvature induced by $\tilde g$ on $\partial M$.

Now, let us construct a function $\tilde R'\in L^{ p}_{\tau-2}(M)$ such that $0\geq \tilde R'\geq R'$, $\tilde R'$ has compact support and if $\tilde Z$ is the zero set of $\tilde R'=0$, $\Yam(\tilde Z, Z_\partial)>0$.

First, consider $E_n:=\displaystyle\bigcup_{i=1}^pf_i(\mathbb{R}^n)\setminus B_n(0)$, where the $f_i$'s are the end charts associated to the asymptotically euclidean manifold $M$. Then, let $\chi_n$ be a smooth bump function characterized by:

i. $\chi_n\equiv 0$ in $\overline{E_{2n}}$,

ii. $\chi_n>0$ in $(M\cup \partial M)\setminus \overline{E_{2n}}$,

iii. $\chi_n\equiv 1$ in $(M\cup \partial M)\setminus \overline{E_{n}}$.

So $\chi_nR'$ has compact support for any $n$, $\chi_nR'\to R'$ pointwise if $n\to\infty$ and if $Z_n$ is the zero set of $\chi_nR'$, $\displaystyle\bigcap_n Z_n=Z$. As a consequence, because of Proposition \ref{2p:continuity_from_above}, $\lambda_\delta(Z_n,Z_\partial )\to \lambda_\delta(Z,Z_\partial)$ for $\delta>\delta^\ast$. Now, the sign of $\lambda_\delta(Z,Z_\partial)$ is the same as the sign of $\Yam(Z, Z_\partial)>0$, hence there is $n_0$ such that $\lambda_\delta(Z_{n_0}, Z_\partial)>0$, and $\Yam(Z_{n_0}, Z_\partial)>0$. Because $0\leq \chi_{n_0}\leq 1$ and $R\leq 0$, just take $\tilde R'=\chi_{n_0}R'$.

As we have both $\tilde R$ and $\tilde R'$ with compact support, using the notation from the previous results in this section we can choose a sequence $\{u_{q_n, r_n}\}_n$ of subcritical minimizers of $F_{q_n,r_n}$ (with $\tilde R$ replacing $R$, $\tilde R'$ replacing $R'$ and $\tilde H$ replacing $H$), $(q_n, r_n)\to (2\bar q, \bar q+1)$, that is uniformly bounded in $W^{2,p}_\tau(M)$. As $W^{2,p}_\tau(M)\hookrightarrow W^{1,2}_{\delta^\ast}(M)$, there is $u\in W^{1,2}_{\delta^\ast}(M)$ such that $u_{q_n, r_n}\to u$ in $W^{1,2}_{\delta^\ast}(M)$ and the convergence is uniform in compact sets. So $u$ is a solution to
\begin{equation}
\label{2eq:critical_equation_tilde}
    \begin{cases}
    -\Delta u+\frac{n-2}{4(n-1)}\tilde R (u+1)=\frac{n-2}{4(n-1)}\tilde R' (u+1)^{2\bar q-1},\ in\ M,\\
\Tr\partial_\nu u+\frac{n-2}{2}\tilde H\Tr (u+1)=\frac{n-2}{2} H'(\Tr (u+1))^{\bar q},\ on\ \partial M,
\end{cases}
\end{equation}
As $\tilde R'\geq R'$, the result follows from Lemma \ref{2l:reduce_curvature}.

Now, assume we can realize $R'$ and $H'$ as scalar curvature in $M$ and mean curvature on $\partial M$, respectively, by a metric $g'$ that is conformal to $g$. The proof that $\Yam(Z,Z_\partial)>0$ is very similar to the case of compact manifolds with boundary (check the proof of Theorem 4.6 in \cite{ST21}). In fact, if $Z$ is bounded or if $Z\setminus K$ has measure zero for some compact set $K$, the proof in the previous paper applies to the asymptotically euclidean case as well, as in that case $W^{1,2}(Z,Z_\partial)=W^{1,2}_\delta(Z,Z_\partial)$.

So assume $Z$ is not bounded and $Z\setminus K$ has positive measure for any compact set $K$. Then if $\Yam(Z,Z_\partial)=0$, $\lambda_\delta(Z,Z_\partial)=0$ and there is $\bar u$ a minimizer for $\lambda_\delta(Z, Z_\partial)$. Just as in the compact case, that implies $\bar u$ is a constant, but the only constant with the proper decay in $W^{1,2}_\delta(Z,Z_\partial)$ in the case that $Z$ is unbounded is $\bar u\equiv0$, and that cannot be true.

\end{proof}
\end{theorem}

This existence result allows us to classify conformal classes of asymptotically euclidean manifolds with boundary with respect to the non-positive curvatures that can be attained in terms of the Yamabe invariant of the manifold. Specifically, we have the following classification.

\begin{theorem}
\label{2t:yamabe_classification}
Let $(M,\partial M, g)$ be an asymptotically euclidean Riemannian manifold with boundary, $g\in W^{2,p}_\tau(M)$, $p>\frac{n}{2}$, $\tau<0$. Then we have that:

i. $\Yam(M, \partial M)>0$ if, and only if, any pair of continuous non-positive functions $(R',H')$, $R'\in L^p_{\tau-2}(M)$, $H'\in W^{1-\frac{1}{p}, p}(\partial M)$, can be attained as the scalar curvature in $M$ and the mean curvature on $\partial M$ induced by a metric $g'$ in the conformal class of $g$.

ii. $\Yam(M, \partial M)=0$ if, and only if, any pair of continuous non-positive functions $(R',H')$, $R'\in L^p_{\tau-2}(M)$, $H'\in W^{1-\frac{1}{p}, p}(\partial M)$, can be attained as the scalar curvature in $M$ and the mean curvature on $\partial M$ induced by a metric $g'$ in the conformal class of $g$, with the exception of the pair $(R'\equiv 0, H'\equiv 0)$.

iii. $\Yam(M, \partial M)<0$ otherwise, that is, if there is a pair of continuous non-positive functions $(R', H')$, $R'\in L^p_{\tau-2}(M)$, $H'\in W^{1-\frac{1}{p}, p}(\partial M)$, with either $R'\not\equiv0$ or $H'\not\equiv 0$ that cannot be attained as the scalar curvature in $M$ and the mean curvature on $\partial M$ by any metric $g'$ on the conformal class of $g$.
\begin{proof}
For simplicity, throughout the proof we will say a pair $(R', H')$ is attainable if there is $g'$ in the conformal class of $g$ such that $R'$ is the scalar curvature in $M$ and $H'$ is the mean curvature on $\partial M$ induced by $g'$. We only have to prove i and ii.

i. Theorem \ref{2t:general_existence_result} implies $\Yam(M,\partial M)>0$ if, and only if, $(R'\equiv 0, H'\equiv 0)$ is attainable, and Lemma \ref{2l:reduce_curvature} implies that $(R'\equiv 0, H'\equiv 0)$ if, and only if, any pair of continuous non-positive functions $(R', H')$ is attainable.

ii. First assume $\Yam(M,\partial M)=0$. By Theorem \ref{2t:general_existence_result}, $(R'\equiv 0, H'\equiv 0)$ is not attainable. Now, let $(R', H')$ be a pair of continuous non-positive functions such that either $R'\not\equiv 0$ or $H'\not\equiv 0$ and $u\not\equiv 0$ be a function supported in $(Z,Z_\partial)$ such that $\lambda_\delta(Z,Z_\partial)=E(u)$ for some $\delta$ (if $W^{1,2}_{\delta^\ast}(Z,Z_\partial)=\{0\}$, $\lambda_\delta(Z,Z_\partial)=+\infty$ and the pair is attainable).

If $E(u)<0$, $\Yam(M,\partial M)<0$, which is not the case. So assume $E(u)=0$. Then $u$ is a minimizer for $\lambda_\delta(M,\partial M)$ as well, and a solution to equation (\ref{2eq:linear_problem}). But $Z\cup Z_\partial \varsubsetneq M\cup\partial M$ by the choice of the pair, hence there is $x\in (M\cup\partial M)\setminus (Z\cup Z_\partial)$ such that $u(x)=0$. However, Lemma 4 of \cite{Max05} would then imply that $u\equiv 0$, contradiction. So actually $\Yam(Z,Z_\partial)>0$ and $(R', H')$ is attainable by Theorem \ref{2t:general_existence_result}.

On the other hand, assume any pair $(R', H')$ of non-negative continuous functions is attainable but for $(R'\equiv 0, H'\equiv 0)$. Then $\Yam(M,\partial M)\leq 0$ by i. If $\Yam(M,\partial M)<0$, let $u$ be a minimizer for $\lambda_\delta(M,\partial M)$ for some $\delta$, and take $\{U_n\}_n$ a sequence of nested bounded open sets such that $\bigcap_n U_n=\emptyset$ and for each $n$, call $V_n$ an open set such that $\bar V_n\varsubsetneq U_n$.

Also, consider $\chi_n$ a smooth non-negative bump function such that $\chi_n\vert_{V_n}\equiv 1$ and $\chi_n\vert_{M\setminus U_n}\equiv 0$. Then, by the dominated convergence theorem, $E((1-\chi_n)u)\to E(u)<0$, so there is $n_0$ such that $E((1-\chi_{n_0})u)<0$. But $v:=(1-\chi_n)u$ vanishes in $V_n$, so $v\in W^{1,2}_{\delta^\ast}(M\setminus V_{n_0}, \partial M)$ and $\Yam(M\setminus V_{n_0}, \partial M)<0$. If we choose $H'\equiv 0$ and $R'$ a bounded, non-positive, smooth function such that $R'<0$ in $V_{n_0}$ and $R'\vert_{M\setminus V_{n_0}}\equiv 0$, Theorem \ref{2t:general_existence_result} tells that $(R', H')$ is NOT attainable, contradicting our hypothesis. So $\Yam(M,\partial M)=0$ instead.
\end{proof}
\end{theorem}

Notice that the proof above shows three facts in addition to the classification result itself that should be highlighted:

i. All proper subsets of a Yamabe positive manifold are Yamabe positive (this also follows from the monotonicity of the invariants by inclusion).

ii. If $(M,\partial M)$ is Yamabe zero, any pair of subsets $(V,V_\partial)$, $V\subset M$, $V_\partial \subset \partial M$ such that either $M\setminus V$ or $\partial M\setminus V_\partial$ contains a relative open set is Yamabe positive,

iii. If $(M,\partial M)$ is Yamabe negative, it is easy to find subsets of $M$ that differ from $M$ by small open sets that are Yamabe negative, thus curvatures that are concentrated on those small open sets are not attainable. It is not hard to see how a similar construction can be used to rule out mean curvatures that are concentrated on small sets in the boundary, or combinations of both.

\section*{Acknowledgements}
We would like to thank David Maxwell for the lead on the proof of Lemma \ref{2l:reduce_curvature}.



\printbibliography

\end{document}